\documentclass[12pt,a4paper,headsepline,abstracton,twoside]{scrartcl}

\usepackage[english]{babel}
\usepackage[utf8]{inputenc}
\usepackage{graphicx}
\usepackage{textcomp}  
\usepackage[T1]{fontenc}

\usepackage[automark]{scrpage2}
\usepackage[a4paper,left=3cm, right=2cm,top=3cm, bottom=3cm]{geometry}

\usepackage{amsthm,amsmath,amssymb,amsfonts}
\usepackage{mathrsfs}

\newtheorem{Theorem}{Theorem}[section]
\newtheorem{Lemma}[Theorem]{Lemma}
\newtheorem{Satz}[Theorem]{Proposition}

\theoremstyle{definition} 
\newtheorem{Bemerkung}[Theorem]{Remark}
\newtheorem{Definition}[Theorem]{Definition}
\newtheorem{Beispiel}[Theorem]{Example}

\newenvironment{bew}{\begin{proof}[Proof]}{\end{proof}}

\begin{document}
\title{Local means and atoms in vector-valued function spaces}
\author{Benjamin Scharf}
\date{October 5, 2010}
\maketitle

	\vspace{-1cm}
\begin{abstract}
The first part of this paper deals with the topic of finding equivalent norms and characterizations for vector-valued Besov and Triebel-Lizorkin spaces $B_{p,q}^s(E)$ and $F_{p,q}^s(E)$. We will deduce general criteria by transferring and extending a theorem of Bui, Paluszy\'{n}ski and Taibleson from the scalar to the vector-valued case. 

By using special norms and characterizations we will derive necessary and sufficient conditions for belonging to a vector-valued function spaces $B_{p,q}^s(E)$ or $F_{p,q}^s(E)$. It will be shown that an element of ${\cal S}'(\mathbb{R}^n,E)$ belongs to a function space if and only if it can be written as a linear combination of harmonic atoms resp.$ $ quarks with suitable conditions for the coefficients. 
\end{abstract}


        \addcontentsline{toc}{section}{Summary}

  	\section{Introduction}
	The aim of this work is to extend the results for atomic and subatomic charaterizations of the function spaces $B_{p,q}^s$ and $F_{p,q}^s$ to the vector-valued function spaces $B_{p,q}^s(E)$ and $F_{p,q}^s(E)$. For a comprehensive treatise of the scalar case ($E=\mathbb{C}$) we refer to chapter 13 and 14 of \cite{Tri97}. A short consideration of the vector-valued case is given in chapter 15 of that book. But the proofs of the crucial theorems 15.8, p.$ $ 114 and 15.11, p.$ $ 116  are only shortly outlined and are mostly based on results for vector-valued function spaces which are well-known in the scalar, but have not yet been considered in the vector-valued case in detail. 

This paper tries to derive these two theorems in wider detail, including the necessary steps before. In chapter 2 we will deal with the fundamentals of vector-valued functions and function spaces. We won't give any proofs mainly because most of them are similar to the scalar case. Many of these were treated in \cite{Tri83}.

In the third chapter we will prove a general result for equivalent norms and characterizations of vector-valued function spaces $B_{p,q}^s(E)$ and $F_{p,q}^s(E)$ in full detail. The scalar version ($E=\mathbb{C}$) of this theorem goes back to Bui, Paluszy\'{n}ski and Taibleson (see \cite{BPT96} and \cite{BPT97}), where the proof, which we will transfer to the vector-valued case, is given in this form in \cite{Ryc99}. Nevertheless, there will be a little modification caused by some minor gap in the original proof. An earlier version with a bit worse, but more general conditions can be found in \cite{Tri92}, section 2.4, p.$ $ 100 for $F_{p,q}^s$ and section 2.5, p.$ $ 132 for $B_{p,q}^s$. In the following we use our result to obtain explicit norms and characterizations which we need to prove atomic and subatomic representations later on.

In the fourth chapter we will derive atomic and subatomic charaterizations for function spaces. We keep close to the approach suggested in \cite{Tri97}, theorem 15.8, p.$ $ 114. Thus we follow chapters 13 and 14 of \cite{Tri97} and transfer the results to the vector-valued case, with minor modifications due to some imperfections in the original proof.
        \addcontentsline{toc}{section}{Einleitung}

	\section{Mathematical fundamentals}
	\subsection{Vector-valued functions and distributions}

Let $E$ be a complex Banach space with norm $\|\cdot|E\|$ and let $E'$ be its dual. With $U_{E}$ we denote the set of all $x \in E$ with $\|x|E\|=1$. Furthermore, let
\begin{align*}
 B_r(x):=\left\{y \in E: \|x-y|E\|\leq r \right\}, \quad B_r:=B_r(0) \quad \text{and} \quad B:=B_1. 
\end{align*}
Let $(M,{\cal M},m)$ be a $\sigma$-finite measure space, which will be the space $\mathbb{R}^n$ with the $\sigma$-algebra of Borel sets and the Lebesgue measure $|\cdot|$ in the sequel. A function $f: M \rightarrow E$ is called $E$-measurable if there exists a subset $M_0$ of $M$ such that $m(M_0)=0$, $f(M \setminus M_0)$ is contained in a separable subspace $E_0$ of $E$ and if the complex-valued functions 
\begin{align*}
 a(f):  x \mapsto a(f(x))
\end{align*}
are measurable for all $a \in E'$. 

If $f$ is $E$-measurable in this sense, then the function $\|f|E\| : M \rightarrow \mathbb{R},\ x \mapsto \|f(x)|E\|$ is measurable because of
\begin{align}
\label{GrundHahn}
 \|f(x)|E\|= \sup_{a \in U_{E'}} |a(f(x))|.
\end{align}
Therefore, we can define the spaces $L_p(E)$ for $0<p\leq \infty$ as follows: 
\begin{align*}
 L_p(M,E):=\left\{f: M \rightarrow E, f \text{ measurable }, \Big\| \|f|E\| \big|L_p(M,{\cal M},m)\Big\|< \infty \right\}.
\end{align*}
We write shortly $L_p(E):=L_p(\mathbb{R}^n,E)$ and $L_p:=L_p(\mathbb{C})$. The spaces $L_p(M,E)$ are (quasi)-Banach spaces. 

For functions $f:M \rightarrow E$ of the form 
\begin{align*}
 f=\sum_{k=1}^K b_k(x) u_k
\end{align*}
with integrable $b_k: M \rightarrow \mathbb{C}$ and $u_k \in E$ for $ k=1,\ldots,K$ we define the Bochner integral as a mapping into $E$ through
\begin{align*}
 \int_{\mathbb{R}^n} f(x) \ dx := \sum_{k=1}^K u_k \int_{\mathbb{R}^n} b_k(x) \ dx.
\end{align*}
For every $a \in E'$ it follows
\begin{align}
\label{GrundLinear}
\begin{split}
 a\left(\int_{\mathbb{R}^n} f(x)\ dx\right)&=a\left(\sum_{k=1}^K u_k \int_{\mathbb{R}^n} b_k(x) \ dx\right)= \sum_{k=1}^K a(u_k) \int_{\mathbb{R}^n} b_k(x) \ dx \\
&= \int_{\mathbb{R}^n}  \sum_{k=1}^K a(u_k) b_k(x) \ dx = \int_{\mathbb{R}^n} a(f(x)) \ dx
\end{split}
\end{align}
and thus with $\eqref{GrundHahn}$
\begin{align}
\label{GrundNormBoch}
 \left\| \int_{\mathbb{R}^n} f(x) \ dx \big|E\right\|\leq \int_{\mathbb{R}^n} \|f(x)|E\| \ dx.
\end{align}
According to that the Bochner integral is a bounded linear operator from the subspace of functions of this form into $E$. This subspace is dense in $L_1(M,E)$ (see $\cite{Gra04}$, section 4.5.c., p.$ $ 318). So the operator can be continued to $L_1(M,E)$ uniquely. We want to call this continuation Bochner integral. Then the properties $\eqref{GrundLinear}$ and $\eqref{GrundNormBoch}$ hold for all $f \in L_1(M,E)$.

We define the Hardy-Littlewood maximal function $M(f)$ for $f \in L_1^{loc}$ as
\begin{align}
 \label{GrundMaximalfunktion}
M(f)(x):=\sup_{B_r(y) \ni x} \frac{1}{|B_r(y)|} \int_{B_r(y)} |f(y)| \ dy.
\end{align}
If, for a given $K:\mathbb{R}^n \rightarrow \mathbb{C}$, there exists a non-negative, monotonically decreasing function $\psi \in L_1((0,\infty))$ with $|K(x)|\leq \psi(|x|)$, then it holds
\begin{align}
 \label{GrundMaxMajo}
\sup_{\delta>0} \left|K_{\delta} * f\right|(x) \leq \|\psi(|\cdot|)|L_1(\mathbb{R}^n,\mathbb{C})\| \cdot M(f)(x)
\end{align}
for $K_{\delta}(x):=\delta^{-n}K(\delta^{-n}x)$ and $f \in L_1^{loc}$. A proof of this proposition can be found in \cite{StW90}, chapter 3, p.$ $ 59. Furthermore, for every $1<p\leq \infty$ there exists a constant $c>0$ such that
\begin{align}
 \label{GrundMaxLp}
\| M(f)|L_p\| \leq c\|f|L_p\|
\end{align}
for all $f \in L_p$ and for every $1<p< \infty$ and $1<q\leq\infty$ there exists a constant $c>0$ such that
\begin{align}
 \label{GrundMaxLplq}
\left\|\left(\sum_{j=1}^{\infty} M(f_j)^q\right)^{\frac{1}{q}}\big|L_p \right\| \leq c 
\left\|\left(\sum_{j=1}^{\infty} |f_j|^q\right)^{\frac{1}{q}}\big|L_p \right\|
\end{align}
for all $\{f_j\}_{j \in \mathbb{N}} \in L_p(l_q)$. References for the proofs are given in \cite{Tri92}, section 2.2.2, p. 89.
\\[1em]
We denote by ${\cal S}(\mathbb{R}^n,E)$ the space of functions $\varphi: \mathbb{R}^n \rightarrow E$ which are infinitely often differentiable and for which the norms 
\begin{align*}
 \|\varphi|E\|_{K,L}:=\sup_{x \in \mathbb{R}^n}  (1+|x|^2)^{\frac{K}{2}} \sum_{|\alpha|\leq L} \|D^{\alpha} \varphi(x)|E\|
\end{align*}
for $K,L \in \mathbb{N}_0$ are finite.
We write shortly ${\cal S}(\mathbb{R}^n):={\cal S}(\mathbb{R}^n,\mathbb{C})$ and, for $\varphi \in {\cal S}(\mathbb{R}^n)$,
\begin{align}
 \label{GrundSNormen}
 \|\varphi\|_{K,L}:=  \|\varphi|\mathbb{C}\|_{K,L}.
\end{align}
The Fourier transform $\hat{\varphi}$ of $\varphi \in {\cal S}(\mathbb{R}^n)$ will be defined as
\begin{align*}
 \hat{\varphi}(\xi):=(2\pi)^{-\frac{n}{2}} \int_{\mathbb{R}^n} \varphi(x) e^{-ix\xi} \ dx,
\end{align*}
whereas we denote the inverse Fourier transform by $\check{\varphi}$. It holds
\begin{align*}
 \check{\varphi}(\xi)=(2\pi)^{-\frac{n}{2}} \int_{\mathbb{R}^n} \varphi(x) e^{ix\xi} \ dx.
\end{align*}
We call a linear map $f: {\cal S}(\mathbb{R}^n) \rightarrow E$ an $E$-valued tempered distribution if there exist constants $c>0$ and $K,L \in \mathbb{N}_0$ such that for all $\varphi \in {\cal S}(\mathbb{R}^n)$ we have
\begin{align*}
 \|f(\varphi)|E\| \leq c \|\varphi\|_{K,L}.
\end{align*}
The set of all this linear maps will be denoted by ${\cal S}'(\mathbb{R}^n,E)$. We say that $f_j$ converges to $f$ in ${\cal S}'(\mathbb{R}^n,E)$ if and only if $f_j(\varphi)$ converges to $f(\varphi)$ for all $\varphi \in {\cal S}(\mathbb{R}^n)$. Such a distribution $f$ will be called regular if there is a measurable, locally Bochner integrable function $g: \mathbb{R}^n \rightarrow E$ so that
\begin{align*}
  f(\varphi)=\int_{\mathbb{R}^n} g(x)\varphi(x) \ dx
\end{align*}
for all $\varphi \in {\cal S}(\mathbb{R}^n)$. As in the scalar case $L_p(E)$ for $1\leq p \leq \infty$ can be understood as a subset of ${\cal S}'(\mathbb{R}^n,E)$ .

For an $f \in {\cal S}'(\mathbb{R}^n,E)$ we define the Fourier transform $\hat{f}$ as
\begin{align*}
 \hat{f} (\varphi):=f(\hat{\varphi}) \text{ for } \varphi \in {\cal S}(\mathbb{R}^n).
\end{align*}
The usual fundamental properties from the scalar case can be transfered.

For $f \in {\cal S}'(\mathbb{R}^n,E)$ and $\psi \in {\cal S}(\mathbb{R}^n)$ we define the convolution as 
\begin{align}
\label{GrundFalt}
 \left(\psi * f\right)(x):=(2\pi)^{-\frac{n}{2}} f\left(\psi(x-\cdot) \right) \text{ for } x \in \mathbb{R}^n,
\end{align}
analogously to the scalar case. The function $\psi*f$ is infinitely often differentiable and there exist $c>0$ and $K,L\in \mathbb{N}_0$ such that
\begin{align}
\label{GrundFalt2}
 \|\left(\psi * f\right)(x)|E\|\leq c(1+|x|^2)^{\frac{K}{2}} \|\psi\|_{K,L}.
\end{align}
As in the scalar case the important relation
\begin{align*}
 \left(\psi*f \right) \hat{\ }= \hat{\psi} \cdot \hat{f}.
\end{align*}
holds.

\subsection{Vector-valued function spaces}
Let $\Omega$ be an open subset of $\mathbb{R}^n$. We set
\begin{align*}
 L_p^{\Omega}(E):=\{f \in L_p(\mathbb{R}^n,E)\cap {\cal S}'(\mathbb{R}^n,E): supp \ \hat{f} \subset \Omega \}
\end{align*}
for $0<p\leq \infty$ and shortly $L_p^{\Omega}$ if $E=\mathbb{C}$. The Nikolskii inequality can be transfered to the vector-valued case, i.e. for $0<p_1<p_2\leq \infty$ there exists a constant $c>0$ such that for all $r>0$ and $f \in L_{p_1}^{B_r}(E)$ it holds
\begin{align}
 \label{GrundNikolskij}
 \|f|L_{p_2}(E)\| \leq c\, r^{n\left(\frac{1}{p_1}-\frac{1}{p_2}\right)} \|f|L_{p_1}(E)\|.
\end{align} 
For a proof see \cite{ScS01}, lemma 1, p. 6. Moreover, let $f \in L_p(E)$ for $1\leq p\leq \infty$ and $g \in L_1$. Then
\begin{align}
 \label{GrundFaltung1}
  \|g*f|L_p(E)\|\leq (2\pi)^{-\frac{n}{2}} \|g|L_1\| \cdot \|f|L_p(E)\|.
\end{align}
If otherwise $0<p<1$, then there exists a constant $c>0$ such that for $f \in L_p^{B_r}(E)$, $g \in L_p^{B_r}$ and $r>0$
\begin{align}
 \label{GrundFaltung2}
  \|g*f|L_p(E)\|\leq c r^{n\left(\frac{1}{p}-1\right)} \|g|L_p\| \cdot \|f|L_p(E)\|
\end{align}
holds. Additionally, one gets for $0<p\leq \infty$, $a>\frac{n}{p}$ and $f \in L_p^{B_r}(E)$
\begin{align}
 \label{GrundPeetre}
\left\| \sup_{z \in \mathbb{R}^n} \frac{\|f(\cdot-z)|E\|}{(1+r|z|)^a}\big|L_p\right\| \leq c \|f|L_p(E)\|.
\end{align}
Proofs for $E=\mathbb{C}$ can be found in \cite{Tri83}, section 1.5.1. resp. 1.4.1.
\\[1em]
Let $\varphi_j$ for $j \in \mathbb{N}_0$ be elements of ${\cal S}(\mathbb{R}^n)$ with
\begin{align}
 \label{GrundResolution}
 \begin{split}
 supp \ \varphi_0 &\subset  \{|\xi| \leq 2\}, \\
 supp \ \varphi_j &\subset\{2^{j-1}\leq |\xi| \leq 2^{j+1}\} \text{ for } j \in \mathbb{N}, \\
 \sum_{j=0}^{\infty} \varphi_j(\xi)&=1 \text{ for all } \xi \in \mathbb{R}^n, \\
 |D^{\alpha} \varphi_j(\xi)| &\leq c_{\alpha} 2^{-j|\alpha|} \text{ for all } \alpha \in \mathbb{N}_0^n.
\end{split}
\end{align}
Then we call $\{\varphi_j\}_{j=0}^{\infty}$ a smooth dyadic resolution of unity.
\begin{Definition}
 \label{GrundDefinitionB}
Let $0<p\leq \infty$, $0<q\leq \infty$, $s \in \mathbb{R}$ and $\{\varphi_j\}_{j=0}^{\infty}$ be a smooth dyadic resolution of unity. For $f\in {\cal S}'(\mathbb{R}^n,E)$ we define
\begin{align*}
 \|f|B_{p,q}^s(E)\|:=\left(\sum_{j=0}^{\infty} 2^{jsq}\|(\varphi_j \hat{f})\check{\ }|L_p(E)\|^q\right)^{\frac{1}{q}}
\end{align*}
(modified if $q=\infty$) and
\begin{align*}
 B_{p,q}^s(E):=\left\{f \in {\cal S}'(\mathbb{R}^n,E): \|f|B_{p,q}^s(E)\|< \infty \right\}.
\end{align*}
\end{Definition}
\begin{Definition}
 \label{GrundDefinitionF}
Let $0<p< \infty$, $0<q\leq \infty$, $s \in \mathbb{R}$ and $\{\varphi_j\}_{j=0}^{\infty}$ be a smooth dyadic resolution of unity. For $f\in {\cal S}'(\mathbb{R}^n,E)$ we define
\begin{align*}
 \|f|F_{p,q}^s(E)\|:=\left\|\left(\sum_{j=0}^{\infty} 2^{jsq} \|(\varphi_j \hat{f})\check{\ }|E\|^q \right)^{\frac{1}{q}}\big|L_p\right\|
\end{align*}
(modified if $q=\infty$) and
\begin{align*}
 F_{p,q}^s(E):=\left\{f \in {\cal S}'(\mathbb{R}^n,E): \|f|F_{p,q}^s(E)\|< \infty \right\}.
\end{align*}
\end{Definition}
We will write shortly $B_{p,q}^s$ for $B_{p,q}^s(\mathbb{C})$ and $F_{p,q}^s$ for $F_{p,q}^s(\mathbb{C})$. As in the scalar case one can show that the definition does not depend on the choice of the smooth dyadic resolution of unity and that the introduced quasi-norms\footnote{In the following we will use the term ``norm`` even if we only have quasi-norms for $p<1$ or $q<1$.} for two different smooth dyadic resolutions of unity are equivalent. Furthermore, the so defined spaces are (quasi)-Banach spaces. We have the fundamental embedding 
\begin{align*}
 B_{p,q_1}^s(E) \hookrightarrow B_{p,q_2}^s(E), \qquad F_{p,q_1}^s(E) \hookrightarrow F_{p,q_2}^s(E)
\end{align*}
for $q_1<q_2$  and
\begin{align*}
 B_{p,\min(p,q)}^s(E) \hookrightarrow F_{p,q}^s(E) \hookrightarrow B_{p,\max(p,q)}^s(E).
\end{align*}
Additionally we have 
\begin{align}
\label{GrundCub}
 B_{\infty,1}^L(E) \hookrightarrow C_{ub}^L(E) \hookrightarrow B_{\infty,\infty}^L(E)
\end{align}
for all $L \in \mathbb{N}_0$, where $C_{ub}^L(E)$ is the set of all $L$ times continuously differentiable functions $f: \mathbb{R}^n \rightarrow E$.

For $m: \mathbb{R}^n \rightarrow \mathbb{C}$ let
\begin{align*}
 \|m\|_N:=\sup_{|\alpha|\leq N} \sup_{x \in \mathbb{R}^n} (1+|x|^2)^{\frac{|\alpha|}{2}} |D^{\alpha} m(x)|.
\end{align*}
Then there exist $c>0$ and $N \in \mathbb{N}$ in dependence of $p$, $q$ and $s$ such that for all infinitely often differentiable functions $m: \mathbb{R}^n \rightarrow \mathbb{C}$ 
\begin{align}
 \label{GrundFourierMult}
  \begin{split}
  \|(m\hat{f})\check{\ }|B_{p,q}^s(E)\| &\leq c\,\|m\|_N \cdot\|f|B_{p,q}^s(E)\| \text{ resp.} \\
  \|(m\hat{f})\check{\ }|F_{p,q}^s(E)\| &\leq c\,\|m\|_N \cdot \|f|F_{p,q}^s(E)\|.
  \end{split}
\end{align}
For a proof in the scalar case see \cite{Tri83}, section 1.5.2., p. 26 and section 1.6.3., p. 31.

Let $I_{\sigma}(f):=((1+|\cdot|^2)^{\frac{\sigma}{2}}\hat{f})\check{\ }$. Then $f$ is an element of $B_{p,q}^s(E)$ if and only if $I_{\sigma}(f)$ is an element of $B_{p,q}^{s-\sigma}(E)$ and we have 
\begin{align}
 \label{GrundLift}
\|\cdot|B_{p,q}^s(E)\| \sim \|I_{\sigma}(\cdot)|B_{p,q}^{s-\sigma}(E)\|,
\end{align}
analogously for $F_{p,q}^s(E)$. A proof in the scalar case can be found in \cite{Tri83}, section 2.3.8., p. 58.

For function spaces the so-called Sobolev embeddings hold: If $0<p_0\leq p_1\leq \infty$, $0<q\leq \infty$, then 
\begin{align}
 \label{GrundSobolevB}
B_{p_0,q}^{s_0} \hookrightarrow B_{p_1,q}^{s_1} \quad \text{if } s_0-\frac{n}{p_0}=s_1-\frac{n}{p_1}.
\end{align}
For a derivation (in the vector-valued case) see e. g. \cite{ScS01}, proposition 3, p. 12. 

If $0<p_0<p_1<\infty$, $0<q_0,q_1\leq \infty$, then
\begin{align}
 \label{GrundSobolevF}
F_{p_0,q_0}^{s_0} \hookrightarrow F_{p_1,q_1}^{s_1} \quad \text{if } s_0-\frac{n}{p_0}=s_1-\frac{n}{p_1}.
\end{align}
The proof for the vector-valued case can be found in \cite{ScS01}, theorem 5, p. 36.

As in the scalar case we define $\mathscr{C}^s(E):=B_{\infty,\infty}^s$ and
\begin{align*}
 \mathscr{C}^{-\infty}(E):=\bigcup_{s \in \mathbb{R}} \mathscr{C}^s(E).
\end{align*}
By $\eqref{GrundSobolevB}$ and $\eqref{GrundSobolevF}$ we have
\begin{align*}
 \mathscr{C}^{-\infty}(E)=\left\{f \in {\cal S}'(\mathbb{R}^n,E): \exists \ p,q,s \text{ with } f \in B_{p,q}^s(E) \vee f \in F_{p,q}^s(E)\right\}.
\end{align*}
Furthermore, we set
\begin{align*}
\sigma_{p}=n\left(\frac{1}{p}-1\right)_+, \hspace{1em} \sigma_{p,q}=n\left(\frac{1}{\min(p,q)}-1\right)_+, 
\end{align*}
where $a_+=\max(a,0)$. Let $\lfloor a \rfloor$ be the biggest integer smaller or equal to $a$ and $\lceil a \rceil$ the smallest integer bigger or equal to $a$.

	\section{Equivalent norms and characterizations for vector-valued
function spaces}
	In the first section of this chapter we will prove a theorem which gives equivalent norms and characterizations for function spaces $B_{p,q}^s(E)$ and $F_{p,q}^s(E)$ in a very general form. In view of notation we stay close to \cite{Tri92} resp. \cite{Tri97} here as well as in the later chapters such that some differences to the proof in \cite{Ryc99}, on which our derivations are based, cannot be avoided.

In the second part we apply the theorem to get explicit equivalent norms and characterizations which we will need later on for our representation by atomic decompositions. 

\subsection{General characterizations}

Let $f: \mathbb{R}^n \rightarrow \mathbb{C}$ be a measurable function. We set $f_j(x):=2^{jn} f(2^{jn}x)$.
\begin{Theorem}
\label{Rychkov}
Let $S+1 \in \mathbb{N}_0$ with
\begin{align}
  \label{Rychkov1} 
S\geq \lfloor s \rfloor,
\end{align}
let $\Psi$, $\psi \in {\cal S}(\mathbb{R}^n)$ and let there be an $\varepsilon>0$ such that
\begin{align}
  \label{Rychkov2} |\Psi(x)|>0 &\text{ for } \left\{|x|<2\varepsilon \right\}, \\
  \label{Rychkov3} |\psi(x)|>0 &\text{ for } \left\{\frac{\varepsilon}{2}<|x|<2\varepsilon \right\}, \\
  \label{Rychkov4} D^{\alpha} \psi(0)=0 &\text{ for } |\alpha|\leq S. 
\end{align}
Furthermore, let $s \in \mathbb{R}$ and
\begin{align}
 \label{RychkovMax}
  \begin{split}(\Psi^*f)_a(x):=&\sup_{y\in \mathbb{R}^n} \frac{\|(\Psi \hat{f})\check{\ }(x-y)|E\|}{(1+|y|)^a}
=\sup_{y\in \mathbb{R}^n} \frac{\|\left(\check{\Psi} * f\right)(x-y)|E\|}{(1+|y|)^a} \\
 (\psi_j^{*}f)_a(x):=&\sup_{y\in \mathbb{R}^n} \frac{\|(\psi(2^{-j}\cdot) \hat{f})\check{\ }(x-y)|E\|}{(1+2^j|y|)^a}
=\sup_{y\in \mathbb{R}^n} \frac{\|\left(\check{\psi}_j * f\right)(x-y)|E\|}{(1+2^j|y|)^a}.
\end{split}
\end{align}

(i) Let $0<p\leq\infty$, $0<q\leq \infty$ and $a>\frac{n}{p}$. Then
\begin{align*}
  \|f|B_{p,q}^s(E)\|_{\Psi,\psi}:= \|(\Psi \hat{f})\check{\ }|L_p(E)\|+\left(\sum_{j=1}^{\infty} 2^{jsq}\|(\psi(2^{-j}\cdot) \hat{f})\check{\ }|L_p(E)\|^q\right)^{\frac{1}{q}}
\end{align*}
and 
\begin{align*}
  \|f|B_{p,q}^s(E)\|_{\Psi,\psi}^a:=\|(\Psi^*f)_a|L_p\|+\left(\sum_{j=1}^{\infty} 2^{jsq}\|(\psi_j^{*}f)_a|L_p\|^q\right)^{\frac{1}{q}}
\end{align*}
(modified in case of $q=\infty$) are equivalent norms for $\|\cdot|B_{p,q}^s(E)\|$. In addition, it holds
\begin{align}
 \label{RychkovCharak1}
B_{p,q}^s(E)=\left\{f \in {\cal S}'(\mathbb{R}^n,E): 
 \|f|B_{p,q}^s(E)\|_{\Psi,\psi}<\infty \right\}
\end{align}
and
\begin{align}
 \label{RychkovCharak2} 
B_{p,q}^s(E)=\left\{f \in {\cal S}'(\mathbb{R}^n,E): 
 \|f|B_{p,q}^s(E)\|_{\Psi,\psi}^a<\infty \right\}.
\end{align}

(ii) Let $0<p<\infty$, $0<q\leq \infty$ and $a>\frac{n}{\min(p,q)}$. Then
\begin{align*}
\|f|F_{p,q}^s(E)\|_{\Psi,\psi}:=\|(\Psi \hat{f})\check{\ }|L_p(E)\|+\left\|\left(\sum_{j=1}^{\infty} 2^{jsq} \|(\psi(2^{-j}\cdot) \hat{f})\check{\ }|E\|^q \right)^{\frac{1}{q}}|L_p\right\|
\end{align*}
and
\begin{align*}
  \|f|F_{p,q}^s(E)\|_{\Psi,\psi}^a:=\|(\Psi^*f)_a|L_p\|+\left\|\left(\sum_{j=1}^{\infty} 2^{jsq}\left((\psi_j^{*}f)_a\right)^q \right)^{\frac{1}{q}}|L_p\right\|
\end{align*}
(modified in case of $q=\infty$) are equivalent norms for $\|\cdot|F_{p,q}^s(E)\|$. In addition, it holds
\begin{align*}
 F_{p,q}^s(E)=\left\{f \in {\cal S}'(\mathbb{R}^n,E): 
 \|f|F_{p,q}^s(E)\|_{\Psi,\psi}<\infty \right\}
\end{align*}
and
\begin{align*}
 F_{p,q}^s(E)=\left\{f \in {\cal S}'(\mathbb{R}^n,E): 
 \|f|F_{p,q}^s(E)\|_{\Psi,\psi}^a<\infty \right\}.
\end{align*}

\end{Theorem}
\begin{bew}
 \bfseries{First step: }\normalfont Let $\Phi$, $\varphi \in {\cal S}(\mathbb{R}^n)$ with
\begin{align}  
\label{RychkovBedPhi}
\begin{split}
 |\Phi(x)|>0 &\text{ for } \left\{|x|<2\varepsilon' \right\}, \\
|\varphi(x)|>0 &\text{ for } \left\{\frac{\varepsilon'}{2}<|x|<2\varepsilon' \right\}
\end{split}
\end{align}
be given and let $(\Phi^*f)_a(x)$ and $(\varphi_j^*f)_a(x)$ be defined analogously as $\eqref{RychkovMax}$. Let $a>0$, $0<p\leq\infty$ ($0<p<\infty$ in case of $F_{p,q}^s(E)$), $0<q\leq \infty$ and $s<S+1$ be fixed. We want to show in this step that there is a constant $C>0$ independent of $f$ such that 
\begin{align}
\label{RychkovFirstPart1}
\begin{split}
  \|(\Psi^*f)_a|L_p\|+\left(\sum_{j=1}^{\infty} 2^{jsq}\|(\psi_j^{*}f)_a|L_p\|^q\right)^{\frac{1}{q}}& \\ \leq C \ 
  \|(\Phi^*f)_a|L_p\|+&C\left(\sum_{j=1}^{\infty} 2^{jsq}\|(\varphi_j^*f)_a|L_p\|^q\right)^{\frac{1}{q}}
\end{split}
\end{align}
and 
\begin{align}
\label{RychkovFirstPart2} 
\begin{split}
\|(\Psi^*f)_a|L_p\|+\left\|\left(\sum_{j=1}^{\infty} 2^{jsq}\left((\psi_j^{*}f)_a\right)^q \right)^{\frac{1}{q}}|L_p\right\|& \\
\leq C \ \|(\Phi^*f)_a|L_p\|+&C\left\|\left(\sum_{j=1}^{\infty} 2^{jsq}(\varphi_j^*f)_a)^q \right)^{\frac{1}{q}}|L_p\right\|
\end{split}
\end{align}
holds. We use the following lemma without a proof here. 
\begin{Lemma}
 \label{ExLambda}
Let $\Phi$, $\varphi \in {\cal S}(\mathbb{R}^n)$ with $\eqref{RychkovBedPhi}$ be given. Then there exist two functions $\Lambda,\lambda \in {\cal S}(\mathbb{R}^n)$ with
\begin{align}
 \label{Rychkovsupport}
 \begin{split}
  supp \ \Lambda \subset& \left\{|x|<2\varepsilon'\right\}, \\
  supp \ \lambda \subset& \left\{\frac{\varepsilon'}{2}<|x|<2\varepsilon'\right\}, 
  \end{split} \\
\label{RychkovLambda3} \Lambda(x)\Phi(x)+&\sum_{j=1}^{\infty} \lambda(2^{-j}x)\varphi(2^{-j}x)=1.
\end{align}
\end{Lemma}
$ $
\\[1em]
For our initial $\Phi$, $\varphi \in {\cal S}(\mathbb{R}^n)$ we choose $\Lambda$, $\lambda \in {\cal S}(\mathbb{R}^n)$ by lemma $\ref{ExLambda}$. Now we multiply $\eqref{RychkovLambda3}$ with $f$, apply the Fourier transform and use properties of the convolution of functions from ${\cal S}(\mathbb{R}^n)$ with elements of ${\cal S}'(\mathbb{R}^n,E)$ (see $\eqref{GrundFalt}$) to get
\begin{align*}
  f= \ \left(\check{\Lambda} * \check{\Phi}\right) * f+ \sum_{k=1}^{\infty}\left(\check{\lambda}_k * \check{\varphi}_k\right) * f
\end{align*}
in ${\cal S}'(\mathbb{R}^n,E)$. Hence we can derive
\begin{align}
\label{RychkovPsi}
   (\check{\psi}_j*f)(y)=\left((\check{\psi}_j*\check{\Lambda}) * (\check{\Phi} * f)\right)(y)+ \sum_{k=1}^{\infty}\left((\check{\psi}_j*\check{\lambda}_k) * (\check{\varphi}_k * f)\right)(y)
\end{align}
for all $y \in \mathbb{R}^n$. With the norm inequality of the Bochner integral (see $\eqref{GrundNormBoch}$) it follows
\begin{align}
\label{RychkovFaltung}
\begin{split}
\|\left(\left(\check{\psi}_j*\check{\lambda}_k\right) * \left(\check{\varphi}_k * f\right)\right) (y)|E\|&\leq \int_{\mathbb{R}^n} \left|\left(\check{\psi}_j*\check{\lambda}_k\right)(z)\right|\cdot\|\left(\check{\varphi}_k * f\right)(y-z)|E\| \ dz \\
&\leq (\varphi_k^*f)_a(y) \int_{\mathbb{R}^n} \left|\left(\check{\psi}_j*\check{\lambda}_k\right)(z)\right| (1+2^k|z|)^a \ dz \\
&\equiv (\varphi_k^*f)_a(y) \cdot I_{j,k}.
\end{split}
\end{align}
The scalar(!) integral $I_{j,k}$ is the same as in \cite{Ryc99}.
\begin{Lemma}
\label{RychkovI}
 Let $\mu, \nu \in {\cal S}(\mathbb{R}^n)$, $M \in \mathbb{Z}, M \geq -1$, $d>0$ and 
\begin{align*}
 D^{\alpha} \mu(0)=0 \text{ for all } \alpha \in \mathbb{N}^n \text{ with } |\alpha|\leq M.
\end{align*}
Then for all $N \in \mathbb{N}$ there exists a constant $C_N$ such that for all $t \in (0,d]$ 
\begin{align*}
 \sup_{z \in \mathbb{R}^n} |\left(\mu(t\cdot)\check{\ } * \check{\nu} \right)(z)|\left(1+|z|\right)^N \leq C_N t^{M+1}.
\end{align*}
\end{Lemma}
\begin{bew} 
A proof can be found in Lemma 1 of \cite{Ryc99}.
\end{bew}
For $k \leq j$ we obtain by the substitution of variables $2^ky \rightarrow y$
\begin{align*}
  \int_{\mathbb{R}^n} \left|\left(\check{\psi}_j*\check{\lambda}_k\right)(z)\right|\cdot (1+2^k|z|)^a \ dz &= \int_{\mathbb{R}^n} 2^{kn}\left|\left(\check{\psi}_{j-k}*\check{\lambda}\right)(2^{k}z)\right|\cdot(1+2^k|z|)^a \ dz \\
  &=\int_{\mathbb{R}^n} \left|\left(\check{\psi}_{j-k}*\check{\lambda}\right)(z)\right| \cdot(1+|z|)^a \ dz \\
  &\leq C_{\psi,\lambda} \sup_{z \in \mathbb{R}^n} \left|\left(\psi(2^{k-j} \cdot)\check{\ }*\check{\lambda}\right)(z)\right| \cdot (1+|z|)^{a+n+1} \\
  &\leq C'_{\psi,\lambda} 2^{(k-j)(S+1)}
\end{align*}
using lemma $\ref{RychkovI}$ with $\mu=\psi$ and $\nu=\lambda$ for $M=S$. In case of $k\geq j$ we deduce 
\begin{align*}
 \int_{\mathbb{R}^n} \left|\left(\check{\psi}_j*\check{\lambda}_k\right)(z)\right|\cdot(1+2^k|z|)^a \ dz &= \int_{\mathbb{R}^n} 2^{jn}\left|\left(\check{\psi}*\check{\lambda}_{k-j}\right)(2^{j}z)\right|\cdot (1+2^k|z|)^a \ dz \\
 &= \int_{\mathbb{R}^n} \left|\left(\check{\psi}*\check{\lambda}_{k-j}\right)(z)\right| \cdot (1+|2^{k-j}z|)^a \ dz \\
&\leq 2^{(k-j)a} \int_{\mathbb{R}^n} \left|\left(\check{\psi}*\check{\lambda}_{k-j}\right)(z)\right|\cdot (1+|z|)^a \ dz \\
&\leq C_{M,\psi,\lambda} 2^{(k-j)a} 2^{(j-k)(M+1)},
\end{align*}
where $M$ can be chosen arbitrarily large since $(D^{\alpha} \lambda) (0)=0$ for all $\alpha \in \mathbb{N}^n$ because of the properties of the support of $\lambda$ (see $\eqref{Rychkovsupport}$). If we choose $M\geq2a-s$, we obtain the estimation
\begin{align}
  \label{RychkovI_j,k}I_{j,k}\leq C_{\lambda,\psi}\left\{
\begin{array}{l l}
 2^{(k-j)(S+1)}&, k \leq j \\
 2^{(j-k)(a-s+1)}&, k\geq j  
\end{array}
\right. .
\end{align}
Furthermore, by definition of the maximal functions in $\eqref{RychkovMax}$
\begin{align*}
 (\varphi_k^*f)_a(y)&\leq (\varphi_k^*f)_a(x) (1+2^k|x-y|)^a \\
&\leq (\varphi_k^*f)_a(x) \max\left(1,2^{(k-j)a}\right)(1+2^j|x-y|)^a.
\end{align*}
If we use this and insert it into $\eqref{RychkovFaltung}$ while applying $\eqref{RychkovI_j,k}$, we get
\begin{align}
\label{RychkovWinner1}
  \sup_{y\in \mathbb{R}^n} \frac{\|\left(\check{\psi}_j*\check{\lambda}_k * \check{\varphi}_k * f\right) (y)|E\|}{(1+2^j|x-y|)^a} &\leq  C_{\psi,\lambda}(\varphi_k^*f)_a(x) \left\{\begin{array}{l l}
 2^{(k-j)(S+1)}&, k \leq j \\
 2^{(j-k)(-s+1)}&, k\geq j  
\end{array}
\right. .
\end{align}
In correspondence, if we replace $\lambda_1$ by $\Lambda$ and $\varphi_1$ by $\Phi$ in the previous calculations, we obtain  
\begin{align}
\label{RychkovWinner2}
   \sup_{y\in \mathbb{R}^n} \frac{\|\left(\check{\psi}_j*\check{\Lambda} * \check{\Phi} * f\right) (y)|E\|}{(1+2^j|x-y|)^a} &\leq C_{\Psi,\Lambda} (\Phi^*f)_a(x) 2^{-j(S+1)}.
\end{align}
One has to keep in mind that only the case $1=k\leq j$ is needed, where we haven't used any conditions of the form $(D^{\alpha} \Lambda) (0)=0$. With the representation of $\check{\psi}_j*f$ in $\eqref{RychkovPsi}$ and with the triangle inequality for $\|\cdot|E\|$ we conclude 
\begin{align*}
  (\psi_j^{*}f)_a(x) \leq C (\Phi^*f)_a(x) \ 2^{-j(S+1)}+C\sum_{k=1}^{\infty} (\varphi_k^*f)_a(x) \left\{\begin{array}{l l}
 2^{(k-j)(S+1)}&,k \leq j \\
 2^{(j-k)(-s+1)}&,k\geq j  
\end{array}\right.  .
\end{align*}
By taking $\delta=\min(S+1-s,1)>0$ (see $\eqref{Rychkov1}$) we arrive at
\begin{align}
  \label{RychkovAbsch1}
  2^{js}(\psi_j^{*}f)_a(x) \leq C 2^{-j\delta} (\Phi^*f)_a(x)+C\sum_{k=1}^{\infty}  2^{ks}(\varphi_k^*f)_a(x) 2^{-|j-k|\delta} .
\end{align}
Analogously, by replacing $\psi_1$ by $\Psi$ in the prior remarks, where we only used the case $k\geq j=1$ and therefore conditions of the form  $\left(D^{\alpha} \Psi\right)(0)=0$ are not necessary, we get
\begin{align}
  \label{RychkovAbsch2}
  (\Psi^*f)_a(x) \leq C (\Phi^*f)_a(x)+C\sum_{k=1}^{\infty}  2^{ks}(\varphi_k^*f)_a(x) 2^{-k\delta}. 
\end{align}
Starting from this pointwise estimates we can now establish our assertions $\eqref{RychkovFirstPart1}$ and $\eqref{RychkovFirstPart2}$. For this we choose a usual method which is applied in \cite{Tri92} several times and which turned into a lemma in \cite{Ryc99}.
\begin{Lemma}
\label{Delta}
  Let $0<p,q\leq \infty$ and $\delta>0$. We assume that for the sequences of $\mathbb{R}$-measurable functions $\{g_k\}_{k=0}^{\infty}$ and $\{G_j\}_{j=0}^{\infty}$ it holds
  \begin{align*}
    |G_j(x)|\leq C_0 \sum_{k=0}^{\infty} 2^{-|k-j|\delta}|g_k(x)| \, \text{ for } \ x \in \mathbb{R}^n,
  \end{align*}
  where $C_0$ is a constant independent of $j$ and $x$. Then there exist constants $C_1$ and $C_2$ (in dependence of $p,q,\delta$) such that
  \begin{align}
    \label{RychkovDelta1}\left(\sum_{j=0}^{\infty} \|G_j|L_p\|^q\right)^{\frac{1}{q}} &\leq C_1 \left(\sum_{j=0}^{\infty} \|g_j|L_p\|^q\right)^{\frac{1}{q}}, \\
   \label{RychkovDelta2} \left\|\left(\sum_{j=0}^{\infty} |G_j|^q\right)^{\frac{1}{q}}|L_p\right\| &\leq C_2 \left\|\left(\sum_{j=0}^{\infty} |g_j|^q\right)^{\frac{1}{q}}|L_p\right\|.
  \end{align}
  \end{Lemma}
\begin{bew}
A proof can be found in Lemma 2 of \cite{Ryc99}. 
\end{bew}
Now we come back to the initial topic. Let $G_0(x):=(\Psi^*f)_a(x)$, $G_j(x)=2^{js}(\psi_j^{*}f)_a(x)$ for $j \in \mathbb{N}$, $g_0(x)=(\Phi^*f)_a(x)$ and $g_k(x)=2^{ks}(\varphi_k^*f)_a(x)$ for $k \in \mathbb{N}$. Then the conditions of lemma $\ref{Delta}$ follow from $\eqref{RychkovAbsch1}$ and $\eqref{RychkovAbsch2}$ and we obtain from $\eqref{RychkovDelta1}$ and $\eqref{RychkovDelta2}$, after slight modification, the desired inequalities $\eqref{RychkovFirstPart1}$ and $\eqref{RychkovFirstPart2}$. 

\bfseries{Second Step:} \normalfont Let $\Psi, \psi \in {\cal S}(\mathbb{R}^n)$ with $\eqref{RychkovBedPhi}$ be given. We want to show that there exists a constant $C>0$ with
\begin{align}
\label{RychkovSchritt2a}
\begin{split}
  \|(\Psi^*f)_a|L_p\|+&\left(\sum_{j=1}^{\infty} 2^{jsq}\|(\psi_j^{*}f)_a|L_p\|^q\right)^{\frac{1}{q}} \\ &\leq C \|(\Psi \hat{f})\check{\ }|L_p(E)\|+C\left(\sum_{j=1}^{\infty} 2^{jsq}\|(\psi(2^{-j}\cdot) \hat{f})\check{\ }|L_p(E)\|^q\right)^{\frac{1}{q}}
\end{split}
\end{align}
and an analogous result for $F_{p,q}^s(E)$. 
At the beginning we choose once again $\Lambda,\lambda \in {\cal S}(\mathbb{R}^n)$ for our given $\Psi, \psi \in {\cal S}(\mathbb{R}^n)$ by lemma $\ref{ExLambda}$ with
\begin{align}
\label{RychkovLambda4}
\begin{split}
 supp \ \Lambda &\subset \left\{|x|<2\varepsilon\right\}, \quad  supp \ \lambda \subset \left\{\frac{\varepsilon}{2}<|x|<2\varepsilon\right\}, \\
 1&= \Lambda(x)\Psi(x)+\sum_{k=1}^{\infty} \lambda(2^{-k}x)\psi(2^{-k}x).
\end{split}
\end{align}
If we replace $x$ by $2^{-j}x$ for $j \in \mathbb{N}$ in the last relation, it follows
\begin{align*}
 1= \Lambda(2^{-j}x)\Psi(2^{-j}x)&+\sum_{k=j+1}^{\infty} \lambda(2^{-k}x)\psi(2^{-k}x)
\end{align*}
and 
\begin{align}
 \label{RychkovDarst}
\left(\check{\psi}_j*f\right)(y)&=\left(\left(\check{\Lambda}_j * \check{\Psi}_j\right) * \left(\check{\psi}_j*f\right)\right)(y)+ \sum_{k=j+1}^{\infty}\left(\left(\check{\psi}_j*\check{\lambda}_k\right) * \left(\check{\psi}_k * f\right)\right)(y)
\end{align}
for all $y \in \mathbb{R}^n$. We deduce for all $N \in \mathbb{N}$ with lemma $\ref{RychkovI}$ ($k\geq j$)
\begin{align*}
 \left|\left(\check{\psi}_j*\check{\lambda}_k\right)(z)\right|&=\left|2^{jn}\left(\check{\psi} * \check{\lambda}_{k-j}\right) (2^{j}z)\right|\leq c_{\psi,\lambda,N}\ \frac{2^{jn}2^{(j-k)N}}{(1+2^{j}|z|)^a}
\end{align*}
(without using any moment conditions on $\psi$) and obviously
\begin{align*}
  \left|\left(\check{\Psi}_j*\check{\Lambda}_j\right)(z)\right|&=2^{jn} \left(\check{\Psi} * \check{\Lambda}\right)(2^jz) \leq c_{\Psi,\Lambda} \frac{2^{jn}}{(1+2^j|z|)^a}.
\end{align*}
If we insert these two estimates into $\eqref{RychkovDarst}$, we obtain
\begin{align}
 \label{RychkovAbschSo}
 \|\left(\check{\psi}_j*f\right)(y)|E\| \leq C_N \sum_{k=j}^{\infty} 2^{jn}2^{(j-k)N} \int_{\mathbb{R}^n} \frac{ \|\left(\check{\psi}_k*f\right)(z)|E\|}{(1+2^j|y-z|)^a} \ dz
\end{align}
for all $f \in {\cal S}'(\mathbb{R}^n,E)$. Now we divide both sides by $(1+2^j|x-y|)^a$ and get
\begin{align*} 
 (\psi_j^{*}f)_a(x) \leq C_N \sum_{k=j}^{\infty} 2^{jn}2^{(j-k)N}  \int_{\mathbb{R}^n} \frac{ \|\left(\check{\psi}_k*f\right)(z)|E\|}{(1+2^j|x-z|)^a} \ dz.
\end{align*}
Let $r \in (0,1]$ be fixed. Keeping in mind $k\geq j$ we arrive with
\begin{align*}
  \frac{\|\left(\check{\psi}_k*f\right)(z)|E\|}{(1+2^j|x-z|)^a}  &\leq \|\left(\check{\psi}_k*f\right)(z)|E\|^r \left[(\psi_k^*f)_a(x)\right]^{1-r} \frac{(1+2^k|x-z|)^{a(1-r)}}{(1+2^j|x-z|)^a} \\
&\leq \|\left(\check{\psi}_k*f\right)(z)|E\|^r \left[(\psi_k^*f)_a(x)\right]^{1-r} \frac{2^{a(k-j)}}{(1+2^k|x-z|)^{ar}}
\end{align*}
 (see $\eqref{RychkovMax}$) at
\begin{align}
\label{RychkovR}
\begin{split}
(\psi_j^{*}f)_a(x) &\leq C_{N} \sum_{k=j}^{\infty} 2^{jn}2^{(j-k)N}2^{a(k-j)}[(\psi_k^*f)_a(x)]^{1-r} \int_{\mathbb{R}^n}  \frac{\|\left(\check{\psi}_k*f\right)(z)|E\|^r}{(1+2^k|x-z|)^{ar}} \ dz \\
&= C_{N} \sum_{k=j}^{\infty} 2^{(j-k)N'}[(\psi_k^*f)_a(x)]^{1-r} \int_{\mathbb{R}^n}  \frac{2^{kn}\|\left(\check{\psi}_k*f\right)(z)|E\|^r}{(1+2^k|x-z|)^{ar}} \ dz,
\end{split}
\end{align}
where $N'=N-a+n$ still can be chosen arbitrarily large. This relation holds in an analogous way for $\Psi$ instead of $\psi_j$ and we get slightly varied
\begin{align*}
 (\Psi^*f)_a(x) \leq C_{N} &[(\Psi^*f)_a(x)]^{1-r} \int_{\mathbb{R}^n} \frac{\|\left(\check{\Psi}*f\right)(z)|E\|^r}{(1+|x-z|)^{ar}} \ dz \\
&+C_{N}\sum_{k=1}^{\infty} 2^{-kN'}[(\psi_k^*f)_a(x)]^{1-r} \int_{\mathbb{R}^n} \frac{2^{kn}\|\left(\check{\psi}_k*f\right)(z)|E\|^r}{(1+2^k|x-z|)^{ar}}\ dz.
\end{align*}
We have to modify these two estimates a bit. For that reason we use a lemma which can be directly adopted from \cite{Ryc99}.
\begin{Lemma}
 \label{RychkovO}
Let $0<r\leq 1$ and $\left\{ b_j\right\}_{j=0}^{\infty}$, $\left\{d_j\right\}_{j=0}^{\infty}$ be two sequences with values in $(0,\infty]$ resp. $(0,\infty)$. Let there be an $N_0 \in \mathbb{N}$ with
\begin{align}
\label{RychkovLemmaO}
 \limsup_{j \rightarrow \infty} \frac{d_j}{2^{jN_0}}< \infty
\end{align}
and for all $N \in \mathbb{N}$ a $C_N>0$ such that 
\begin{align*}
 d_j \leq C_N \sum_{k=j}^{\infty} 2^{(j-k)N}b_kd_k^{1-r} \hspace{ 2em } \text{ for } j \in \mathbb{N}_0
\end{align*}
holds. Then for all $N \in \mathbb{N}$ we have
\begin{align*}
 d_j^r \leq C_N \sum_{k=j}^{\infty} 2^{(j-k)Nr}b_k \hspace{ 2em }\text{ for } j \in \mathbb{N}_0
\end{align*}
with the same constants $C_N$.
\end{Lemma}
\begin{bew}
A proof is given in lemma 3 of \cite{Ryc99}.
\end{bew}

For fixed $x \in \mathbb{R}^n$ we make use of lemma $\ref{RychkovO}$ with $d_j=(\psi_j^*f)_a(x)$ for $j \in \mathbb{N}$, $d_0=(\Psi^*f)_a(x)$ and
\begin{align*}
 b_j=\int_{\mathbb{R}^n} \frac{2^{kn}\|\left(\check{\psi}_k*f\right)(z)|E\|^r}{(1+2^k|x-z|)^{ar}}\ dz \ \text{ for all } j \in \mathbb{N}, \hspace{1 em} b_0=\int_{\mathbb{R}^n} \frac{\|\left(\check{\Psi}*f\right)(z)|E\|^r}{(1+|x-z|)^{ar}} \ dz.
\end{align*}
We want to point out that we vary the procedure from $\cite{Ryc99}$ a bit here. We deal with the question whether the $d_j$ fulfil condition $\eqref{RychkovLemmaO}$ in the last step of the proof. The other conditions precisely result from the calculations above (see $\eqref{RychkovR}$). If applicable, we get
\begin{align}
\label{RychkovAbsch}   
(\psi_j^{*}f)_a(x)^r \leq C_{N}' \sum_{k=j}^{\infty} 2^{(j-k)Nr} \int_{\mathbb{R}^n}  \frac{2^{kn}\|\left(\check{\psi}_k*f\right)(z)|E\|^r}{(1+2^k|x-z|)^{ar}} \ dz
\end{align}
and 
\begin{align}
\label{RychkovAbsch4}
 \begin{split} (\Psi^*f)_a(x)^r \leq C_{N}' \int_{\mathbb{R}^n} &\frac{\|\left(\check{\Psi}*f\right)(z)|E\|^r}{(1+|x-z|)^{ar}} \ dz \\
   &+C_{N}'\sum_{k=1}^{\infty} 2^{-kNr} \int_{\mathbb{R}^n} \frac{2^{kn}\|\left(\check{\psi}_k*f\right)(z)|E\|^r}{(1+2^k|x-z|)^{ar}}\ dz
\end{split}
\end{align}
with $C_{N}'= C_{N+a-n}$. Here  $C_{N}'$ does not depend on $f \in {\cal S}'(\mathbb{R}^n,E)$, $j \in \mathbb{N}$ or $r \in (0,1]$. 

We like to note that $\eqref{RychkovAbsch}$ in the case $r>1$ can more easily be derived from $\eqref{RychkovAbschSo}$ if we replace $a$ by $a+n+1$. By applying Hölder's inequality two times we arrive at
\begin{align*}
  \|&\left(\check{\psi}_j*f\right)(y)|E\|  \\
  &\leq C_N \sum_{k=j}^{\infty} 2^{jn}2^{(j-k)N} \int_{\mathbb{R}^n} \frac{ \|\left(\check{\psi}_k*f\right)(z)|E\|}{(1+2^j|y-z|)^{a+n+1}} \ dz \\
&\leq C_N' \left(\sum_{k=j}^{\infty} 2^{(j-k)(N-1+\frac{n}{r}-a)r} \int_{\mathbb{R}^n} \frac{2^{kn} \|\left(\check{\psi}_k*f\right)(z)|E\|^r}{(1+2^k|y-z|)^{ar}} \ dz \right)^{\frac{1}{r}}.
\end{align*}
If we use the inequality
\begin{align*}
  (1+2^j|x-z|)^a \leq (1+2^j|x-y|)^a(1+2^j|y-z|)^a
\end{align*}
when dividing by $(1+2^j|x-y|)^a$, we get
\begin{align*}
  (\psi_j^{*}f)_a(x) \leq C_N' \left(\sum_{k=j}^{\infty} 2^{(j-k)(N-1+\frac{n}{r}-a)r} \int_{\mathbb{R}^n} \frac{2^{kn} \|\left(\check{\psi}_k*f\right)(z)|E\|^r}{(1+2^k|x-z|)^{ar}} \ dz \right)^{\frac{1}{r}}
\end{align*}
and an analogous result for $ (\Psi^*f)_a(x)$, which provides the desired results $\eqref{RychkovAbsch}$ and $\eqref{RychkovAbsch4}$ - because $N \in \mathbb{N}$ was arbitrary - in case of $r>1$.

By our assumptions on $a$ we can choose $r$ in such a way that $\frac{n}{a}<r<p$ resp. $\frac{n}{a}<r<\min(p,q)$. Then we have $h(x):=\frac{1}{(1+|x|)^{ar}} \in L_1$. The majority property (see $\eqref{GrundMaxMajo}$) yields for all $t>0$ 
\begin{align*}
 |\left(h_t * g\right)(x)| \leq c\, M(g)(x)
\end{align*}
for the Hardy-Littlewood maximal function $M(g)$ introduced in $\eqref{GrundMaximalfunktion}$. If we use this for $\eqref{RychkovAbsch}$ (and $\eqref{RychkovAbsch4}$) with $g(z)= \|\left(\check{\psi}_k*f\right)(z)|E\|^r$ and $N=\lfloor \max(-s,0) \rfloor+1$, we come to
\begin{align*}
2^{jsr}(\psi_j^{*}f)_a(x)^r \leq C_{N}' \sum_{k=j}^{\infty} 2^{(j-k)\delta} 2^{ksr}M\left(\|\left(\check{\psi}_k*f\right)|E\|^r\right)(x)
\end{align*}
and to an analogous result for $(\Psi^*f)_a(x)^r$ with a suitable $\delta>0$. Now we apply lemma $\ref{Delta}$ with $G_j=2^{jsr}[(\psi_j^{*}f)_a]^r$ for $j \in \mathbb{N}$, $G_0=[(\Psi^*f)_a] ^r$, $g_k=2^{ksr}M\left(\|\left(\check{\psi}_k*f\right)(z)\|^r\right)$ for $k \in \mathbb{N}$, $g_0=M\left(\|\left(\check{\Psi}*f\right)(z)\|^r\right)$, $\tilde{q}=\frac{q}{r}$ and $\tilde{p}=\frac{p}{r}$. We obtain
\begin{align*}
  \|(\Psi^*f)_a|L_p\|+&\left(\sum_{j=1}^{\infty} 2^{jsq}\|(\psi_j^{*}f)_a|L_p\|^q\right)^{\frac{1}{q}} \leq C \left\|M\left(\|(\Psi \hat{f})\check{\ }|E\|^r\right)\big|L_{\frac{p}{r}}\right\|^{\frac{1}{r}}\\
&\hspace{5em}+C\left(\left(\sum_{j=1}^{\infty} 2^{jsq}\left\|M\left(\|(\psi(2^{-j}\cdot) \hat{f})\check{\ }|E\|^r\right)\big|L_{\frac{p}{r}}\right\|^{\frac{q}{r}}\right)^{\frac{r}{q}}\right)^{\frac{1}{r}}
\end{align*}
resp.$ $ the $F_{p,q}^s(E)$-analogue. Because $\frac{p}{r}>1$ and in case of $F_{p,q}^s(E)$ as well $\frac{q}{r}>1$ (and $\frac{p}{r}<\infty$) it follows from the boundedness of the maximal operator from $L_{\frac{p}{r}}$ to $L_{\frac{p}{r}}$ resp. from $L_{\frac{p}{r}}(l_{\frac{q}{r}})$ to $L_{\frac{p}{r}}(l_{\frac{q}{r}})$ (see $\eqref{GrundMaxLp}$ resp.$ $ $\eqref{GrundMaxLplq}$)
\begin{align*}
  \|(\Psi^*f)_a|L_p\|+&\left(\sum_{j=1}^{\infty} 2^{jsq}\|(\psi_j^{*}f)_a|L_p\|^q\right)^{\frac{1}{q}}  \leq C' \left\|\ \|(\Psi \hat{f})\check{\ }|E\|^r\ |L_{\frac{p}{r}}\right\|^{\frac{1}{r}}\\
&\hspace{6em}+C'\left(\left(\sum_{j=1}^{\infty} 2^{jsq}\left\|\ \|(\psi(2^{-j}\cdot) \hat{f})\check{\ }|E\|^r\ \big|L_{\frac{p}{r}}\right\|^{\frac{q}{r}}\right)^{\frac{r}{q}}\right)^{\frac{1}{r}}
\end{align*}
resp.$ $ the $F_{p,q}^s(E)$-analogue which matches (at close view) our desired result $\eqref{RychkovSchritt2a}$.

\bfseries{Third step:} \normalfont Now we will conclude the equivalences of the norms $\|\cdot|B_{p,q}^s(E)\|$, $\|\cdot|B_{p,q}^s(E)\|_{\Psi,\psi}$ and $\|\cdot|B_{p,q}^s(E)\|_{\Psi,\psi}^*$ by the results of the first and the second step. We choose a smooth dyadic resolution of unity consisting of the functions $\Phi=\varphi_0$ and $\{\varphi_j\}_{j \in \mathbb{N}}$ with $\varphi_j(\cdot)=\varphi(2^{-j}\cdot)$ (see $\eqref{GrundResolution}$) with 
\begin{align*}
 \varphi_0(\xi)>0 \text{ for } |\xi|<2 , \hspace{1em} \varphi(\xi)>0 \text{ for } \frac{1}{2}<|\xi|<2.
\end{align*}
Obviously it holds
\begin{align*}
\|(\Psi \hat{f})\check{\ }|L_p(E)\|&+\left(\sum_{j=1}^{\infty} 2^{jsq}\|(\psi(2^{-j}\cdot) \hat{f})\check{\ }|L_p(E)\|^q\right)^{\frac{1}{p}} 
\\ &\hspace{8em}\leq 
  \|(\Psi^*f)_a|L_p\|+\left(\sum_{j=1}^{\infty} 2^{jsq}\|(\psi_j^{*}f)_a|L_p\|^q\right)^{\frac{1}{p}}.
\end{align*}
By the first step (see $\eqref{RychkovFirstPart1}$) we get
\begin{align*}
  \|(\Psi^*f)_a|L_p\|+&\left(\sum_{j=1}^{\infty} 2^{jsq}\|(\psi_j^{*}f)_a|L_p\|^q\right)^{\frac{1}{p}} \\ &\hspace{6 em} \leq C \|(\Phi^*f)_a|L_p\|+C\left(\sum_{j=1}^{\infty} 2^{jsq}\|(\varphi_j^*f)_a|L_p\|^q\right)^{\frac{1}{p}}
\end{align*}
because $\Phi$ and $\varphi$ fulfil the necessary conditions for $\varepsilon'=1$. Now it follows by the second step, applied to $\Phi$ and $\varphi$,
\begin{align*}
\|(\Phi^*f)_a|L_p\|&+\left(\sum_{j=1}^{\infty} 2^{jsq}\|(\varphi_j^*f)_a|L_p\|^q\right)^{\frac{1}{p}}\\ 
&\leq C' \|(\Phi \hat{f})\check{\ }|L_p(E)\|+C'\left(\sum_{j=1}^{\infty} 2^{jsq}\|(\varphi(2^{-j}\cdot) \hat{f})\check{\ }|L_p(E)\|^q\right)^{\frac{1}{q}}= C'\|f|B_{p,q}^s(E)\|
\end{align*}
if our not yet proven condition of finiteness on $d_j$ in lemma $\ref{RychkovO}$ is true. We will turn our attention to this question immediately.

Otherwise we obtain from the first step - this time by interchanging the roles of $\varphi$ and $\psi$ resp. $\Phi$ and $\Psi$ (this can be done, because $D^{\alpha} \varphi(0)=0$ for all $\alpha \in \mathbb{N}^n$) - and from the second step, applied to $\Psi$ and $\psi$,
\begin{align*}
  \|f|B_{p,q}^s(E)\| &\leq \|(\Phi^*f)_a|L_p\|+\left(\sum_{j=1}^{\infty} 2^{jsq}\|(\varphi_j^*f)_a|L_p\|^q\right)^{\frac{1}{p}} \\
&\leq C\|(\Psi^*f)_a|L_p\|+C\left(\sum_{j=1}^{\infty} 2^{jsq}\|(\psi_j^{*}f)_a|L_p\|^q\right)^{\frac{1}{p}} \\ 
&\leq C'\|(\Psi \hat{f})\check{\ }|L_p(E)\|+C'\left(\sum_{j=1}^{\infty} 2^{jsq}\|(\psi(2^{-j}\cdot) \hat{f})\check{\ }|L_p(E)\|^q\right)^{\frac{1}{q}}
\end{align*}
if our not yet proven condition on the finiteness of $d_j$ in lemma $\ref{RychkovO}$ is true. 

Now let's take a look at this condition: Let us at first assume $f$ to be in $B_{p,q}^s(E)$. Then by the lift property (see $\eqref{GrundLift}$) and the Sobolev embeddings (see $\eqref{GrundSobolevB}$) there is a $\sigma \in \mathbb{N}$ such that  $g:=((1+|\xi|^2)^{-\sigma}\hat{f})\check{\ } \in L_{\infty}(E)$. We obtain
\begin{align*}
 d_j=(\psi_j^*f)_a(x)&=\sup_{y \in \mathbb{R}^n} \frac{\left\|\left(\psi(2^{-j} \cdot) \hat{f}\right)\check{ } \ (y)\big|E\right\|}{(1+2^{j}|x-y|)^a} \\
		     &=\sup_{y \in \mathbb{R}^n} \frac{\left\|\left(\psi(2^{-j}\cdot) (1+|\xi|^2)^{+\sigma} (1+|\xi|^2)^{-\sigma}\hat{f}\right)^{\vee}(y)\big|E\right\|}{(1+2^{j}|x-y|)^a} \\
    &= \sup_{y \in \mathbb{R}^n} \frac{\left\|\left((\psi(2^{-j}\cdot) (1+|\xi|^2)^{\sigma})\check{\ } *g\right)(y)\big|E\right\|}{(1+2^{j}|x-y|)^a} \\
\intertext{By properties of the convolution we estimate by}
    d_j=(\psi_j^*f)_a(x)&\leq \|g|L_{\infty}(E)\|\cdot \|(\psi(2^{-j}\cdot) (1+|\xi|^2)^{\sigma})\check{\ }|L_1\| \\
    &\leq C \|g|L_{\infty}(E)\| \cdot \sum_{|\alpha| \leq 2\sigma} \|D^{\alpha} \check{\psi}_j|L_1\| \\
    &\leq C'2^{j\sigma} \|g|L_{\infty}(E)\|  \cdot \sum_{|\alpha| \leq 2\sigma}\|D^{\alpha} \check{\psi}|L_1\|,
\end{align*}
Therefore, we get the requested condition $\eqref{RychkovLemmaO}$ with $N_0=2\sigma$ and the desired result follows. The proof in the $F_{p,q}^s(E)$-case is the same.

\bfseries{Fourth step:} \normalfont Last but not least we show the characterizations $\eqref{RychkovCharak1}$ and $\eqref{RychkovCharak2}$ for $B_{p,q}^s(E)$. The proof of the $F_{p,q}^s(E)$-case is the same.

In the first step we didn't use the condition $f \in B_{p,q}^s(E)$. So, if for an $f \in {\cal}S'(\mathbb{R}^n,E)$ we have $\|f|B_{p,q}^s(E)\|< \infty$, then $\|f|B_{p,q}^s(E)\|_{\psi,\Psi}^{a}<\infty$ for admissible $a$ and vice versa. Therefore we have $\eqref{RychkovCharak2}$.

In the second step we used the condition $f \in B_{p,q}^s(E)$ only in lemma $\ref{RychkovO}$ for fulfilling $\eqref{RychkovLemmaO}$. If we just assume $f \in {\cal S}'(\mathbb{R}^n,E)$ instead, then there exist constants $c>0, K,L \in \mathbb{N}_0$ such that for all $\varphi \in {\cal S}(\mathbb{R}^n)$ it holds 
\begin{align*}
  \|\left(\check{\varphi}*f\right)(x)|E\|\leq c\ (1+|x|^2)^{\frac{K}{2}} \|\check{\varphi}\|_{K,L},
\end{align*}
see $\eqref{GrundFalt2}$. Hence it follows if $a\geq K$ 
\begin{align*}
 d_j=(\psi_j^*f)_a(x)&=\sup_{y \in \mathbb{R}^n} \frac{\|\left(\check{\psi}_j*f\right)(x-y)|E\|}{(1+2^j|y|)^a} \\
&\leq \sup_{y \in \mathbb{R}^n} \frac{(1+|x-y|^2)^{\frac{K}{2}}}{(1+2^j|y|)^a} \|\check{\psi}_j\|_{K,L} \\
&\leq c'(1+|x|^2)^{\frac{K}{2}} 2^{j(L+n)} \|\check{\psi}\|_{K,L},
\end{align*}
where $c'$ is independent of $j$ and $x$. So the conditions of lemma $\ref{RychkovO}$ are fulfilled for ``large`` $a$ with $N_0=L+n$. Thus it follows that if $\|f|B_{p,q}^s(E)\|_{\psi,\Psi}$ is finite, then as well $\|f|B_{p,q}^s(E)\|_{\psi,\Psi}^{a}$ is finite for these $a$ and hence also $f \in B_{p,q}^s(E)$ by the third step. If $\|f|B_{p,q}^s(E)\|<\infty$, then  $\|f|B_{p,q}^s(E)\|_{\psi,\Psi}^{a}$ is finite for admissible $a$ by the first step and the end of the third step and therefore obviously $\|f|B_{p,q}^s(E)\|_{\psi,\Psi}$, too.
\end{bew}

From the above theorem and its proof a proposition on \textbf{continuous} versions of the norm follows by slight, but technically complex modifications. An example for such norms is given in \cite{Tri92}, section 2.4.1, p.$ $ 101 and section 2.4.3, p.$ $ 115. 
\begin{Satz}
\label{RychkovFolgerung}
Let $d>0$.

(i) Under the assumptions of theorem $\ref{Rychkov}$, part (i) for $\Psi$ and $\psi$
\begin{align*}
  \|f|B_{p,q}^s(E)\|_{\Psi,\psi}^I:=\|(\Psi \hat{f})\check{\ }|L_p(E)\|+\left(\int_0^{1} t^{-sq}\|(\psi(t\cdot) \hat{f})\check{\ }|L_p(E)\|^q \ \frac{dt}{t} \right)^{\frac{1}{q}}
\end{align*}
and 
\begin{align*}
  \|f|B_{p,q}^s(E)\|_{\Psi,\psi}^{\sup}:=\|\sup_{|\cdot-y|\leq d}\|&(\Psi \hat{f})\check{\ }(\cdot)|E\| \ |L_p\| \\
&+\left(\int_0^{1} t^{-sq}\|\sup \| (\psi(\tau\cdot) \hat{f})\check{\ }|E\| \ |L_p\|^q \ \frac{dt}{t}\right)^{\frac{1}{q}}
\end{align*}
are equivalent norms for $\|\cdot|B_{p,q}^s(E)\|$, where $\sup$ is the supremum taken over \\ $\left\{|x-y|\leq dt,t\leq \tau \leq 2t\right\}$ for a fixed $x \in \mathbb{R}^n$. It holds
\begin{align*}
B_{p,q}^s(E)=\left\{f \in {\cal S}'(\mathbb{R}^n,E): 
 \|f|B_{p,q}^s(E)\|_{\Psi,\psi}^I<\infty \right\}
\end{align*}
and
\begin{align*}
B_{p,q}^s(E)=\left\{f \in {\cal S}'(\mathbb{R}^n,E): 
 \|f|B_{p,q}^s(E)\|_{\Psi,\psi}^{\sup}<\infty \right\}.
\end{align*}

(ii) Under the assumptions of theorem $\ref{Rychkov}$, part (ii) for $\Psi$ and $\psi$
\begin{align}
  \label{RychkovFolg1}
 \|f|F_{p,q}^s(E)\|_{\Psi,\psi}^{I} \|(\Psi \hat{f})\check{\ }|L_p(E)\|+\left\|\left(\int_0^1 t^{-sq} \|(\psi(t\cdot) \hat{f})\check{\ }|E\|^q \ \frac{dt}{t} \right)^{\frac{1}{q}}|L_p\right\|
\end{align}
and
\begin{align}
  \label{RychkovFolg2}
 \begin{split} \|f|F_{p,q}^s(E)\|_{\Psi,\psi}^{\sup}:=\|\sup_{|\cdot-y|\leq d}\|(&\Psi \hat{f})\check{\ }(\cdot)|E\| \ |L_p\| \\
&+\left\|\left(\int_0^1 t^{-sq}\sup \|(\psi(\tau\cdot) \hat{f})\check{\ }|E\|^q \right)^{\frac{1}{q}}|L_p\right\|
\end{split}
\end{align}
are equivalent norms for $\|\cdot|B_{p,q}^s(E)\|$, where $\sup$ is the supremum taken over \\ $\left\{|x-y|\leq dt,t\leq \tau \leq 2t\right\}$ for a fixed $x \in \mathbb{R}^n$. It holds
\begin{align}
 \label{RychkovFolgCharak1}
F_{p,q}^s(E)=\left\{f \in {\cal S}'(\mathbb{R}^n,E): 
 \|f|F_{p,q}^s(E)\|_{\Psi,\psi}^I<\infty \right\}
\end{align}
and
\begin{align}
 \label{RychkovFolgCharak2}
F_{p,q}^s(E)=\left\{f \in {\cal S}'(\mathbb{R}^n,E): 
 \|f|F_{p,q}^s(E)\|_{\Psi,\psi}^{\sup}<\infty \right\}.
\end{align}

\end{Satz}
\begin{bew}
We restrict ourselves to the case of the $F_{p,q}^s(E)$-spaces, the $B_{p,q}^s(E)$-case can be treated in an analogous way. For this purpose we first consider $\eqref{RychkovFolg2}$, which is obviously larger than $\eqref{RychkovFolg1}$, and show that we can estimate this term by $C\|f|F_{p,q}^s(E)\|$. On that account we look back at the first step of the proof of theorem \ref{Rychkov}. But this time we start with $\psi (\tau \cdot)$ with $1 \leq \tau \leq 4$ instead of $\psi$. For given $\Phi$ and $\varphi$ we again choose associated $\Lambda$ and $\lambda$ by lemma $\ref{ExLambda}$. We argue in the same way as in $\eqref{RychkovPsi}$ and $\eqref{RychkovFaltung}$ and obtain
  \begin{align*}
   \|\left(\psi(\tau \cdot)\check{\ }_j*\check{\lambda}_k * \check{\varphi}_k * f\right)(y)|E\|\leq (\varphi_k^*f)_a(y) \int_{\mathbb{R}^n} \left|\left(\psi(\tau \cdot)\check{\ }_j*\check{\lambda}_k\right)(z)\right|(1+2^k|z|)^a \ dz. 
  \end{align*}
  Thereby observe $\psi(\tau \cdot)\check{\ }_j=\psi(\tau 2^{-j}\cdot)\check{\ }$. Now we apply lemma $\ref{RychkovI}$ as in step 1. For $k \leq j$ we get  with the substitution $2^kz \rightarrow z$
  \begin{align*}
    \int \limits_{\mathbb{R}^n} \left|\left(\psi(\tau \cdot)\check{\ }_j*\check{\lambda}_k\right)(z)\right|(1+2^k|z|)^a \ dz &= \int \limits_{\mathbb{R}^n} 2^{kn}\left|\left(\psi(\tau \cdot)\check{\ }_{j-k}*\check{\lambda}\right)(2^{k}z)\right|(1+2^k|z|)^a \ dz \\
  &\leq c \sup_{z \in \mathbb{R}^n} \left|\left(\psi(2^{k-j}\tau \cdot)\check{\ }*\check{\lambda}\right)(z)\right| \cdot (1+|z|)^{a+n+1} \\
  &\leq c' 2^{(k-j)(S+1)}.
  \end{align*}
  In case of $k \geq j$ we obtain by the substitution $2^j\tau^{-1}y \rightarrow y$ and by an analogous calculation as in the proof of the theorem
  \begin{align*}
 \int_{\mathbb{R}^n} &\left|\left(\psi(\tau \cdot)\check{\ }_j*\check{\lambda}_k\right)(z)\right| \cdot (1+2^k|z|)^a \ dz \\
&\hspace{5em}= \int_{\mathbb{R}^n} \tau^{-n}2^{jn}\left|\left(\check{\psi}*\lambda(2^{j-k}\tau^{-1}\cdot)\check{\ }\right)(2^{j}\tau^{-1}z)\right|\cdot (1+2^k|z|)^a \ dz \\
 &\hspace{5em}= \int_{\mathbb{R}^n} \left|\left(\check{\psi}*\lambda(2^{j-k}\tau^{-1}\cdot)\check{\ }\right)(z)\right|\cdot (1+2^{k-j}\tau |z|)^a \ dz \\
& \hspace{5em}\leq c_M 2^{(k-j)a} 2^{(j-k)(M+1)},
\end{align*}
where $c_M$ and $c'$ do not depend on $\tau$. Hence this results in a counterpart of $\eqref{RychkovWinner1}$
\begin{align*}
 \sup_{y\in \mathbb{R}^n} \frac{\|\left(\psi(\tau \cdot)\check{\ }_j*\check{\lambda}_k * \check{\varphi}_k * f \right)(y)|E\|}{(1+2^j|x-y|)^a} &\leq  C_{\psi,\lambda}\ \varphi_{k}^*(x) \left\{\begin{array}{l l}
 2^{(k-j)(S+1)}&, k \leq j \\
 2^{(j-k)(-s+1)}&, k\geq j  \end{array}
\right. 
\end{align*}
independent of $\tau \in [1,4]$. We again come to 
\begin{align*}
 2^{js} \sup_{y\in \mathbb{R}^n} \frac{\|\left(\psi(\tau \cdot)\check{\ }_j* f\right) (y)|E\|}{(1+2^j|x-y|)^a} \leq  C 2^{-j\delta} (\Phi^*f)_a(x)+C\sum_{k=1}^{\infty}  2^{ks}(\varphi_k^*f)_a(x) 2^{-|j-k|\delta}  
\end{align*}
with $\delta= \min(S+1-s,1)>0$ and (taken over from step 1)
\begin{align*}
 \sup_{y\in \mathbb{R}^n} \frac{\|\left(\check{\Psi}_j* f \right)(y)|E\|}{(1+|x-y|)^a} \leq C (\Phi^*f)_a(x)+C\sum_{k=1}^{\infty}  2^{ks}(\varphi_k^*f)_a(x) 2^{-k\delta}. 
\end{align*}
If we restrict each supremum to the domain $|x-y|\leq d2^{-j+1}$ and use that for these $y$ the inequality $(1+2^j|x-y|)^a \leq c_d$  with a constant $c_d>0$ independent of $j$ holds, we get
\begin{align*}
 2^{js} \sup_{\underset{1\leq \tau \leq 4}{|x-y|\leq d2^{-j+1},}}\!\! \|\left(\psi(\tau \cdot)\check{\ }_j* f\right) (y)|E\| \leq  \! C' 2^{-j\delta} (\Phi^*f)_a(x)+\!C'\sum_{k=1}^{\infty}  2^{ks}(\varphi_k^*f)_a(x) 2^{-|j-k|\delta}  
\end{align*}
and
\begin{align*}
 \sup_{|x-y|\leq d} \|\left(\check{\Psi}_j* f\right)(y)|E\| \leq C' (\Phi^*f)_a(x)+C'\sum_{k=1}^{\infty}  2^{ks}(\varphi_k^*f)_a(x) 2^{-k\delta}. 
\end{align*}
By applying lemma $\ref{Delta}$ as in step 1 this yields
\begin{align*}
  \left\|\sup_{|\cdot-y|\leq d}\|(\Psi \hat{f})\check{\ }(\cdot)|E\| \ \big|L_p\right\|&+\left\|\left(\sum_{j=1}^{\infty} 2^{jsq}\sup_{\underset{1\leq \tau \leq 4}{|x-y| \leq d2^{-j+1},}} \|(\psi(2^{-j}\tau\cdot) \hat{f})\check{\ }|E\|^q \right)^{\frac{1}{q}}\big|L_p\right\| \\ 
\leq c \|(\Phi^*f&)_a|L_p\|+c\left\|\left(\sum_{j=1}^{\infty} 2^{jsq}(\varphi_j^*f)_a)^q \right)^{\frac{1}{q}}\big|L_p\right\|.
\end{align*}
But now we have for all $j \in \mathbb{N}$
\begin{align*}
 \int_{2^{-j}}^{2^{-j+1}} t^{-sq} \sup \|(\psi(t\cdot) \hat{f})\check{\ }|E\|^q \ \frac{dt}{t} \leq c_0 2^{jsq}\sup_{\underset{1\leq \tau \leq 4}{|x-y|\leq d2^{-j+1},}} \|(\psi(2^{-j}\tau\cdot) \hat{f})\check{\ }|E\|^q,
\end{align*}
where $\sup$ is the supremum for a fixed $x \in \mathbb{R}^n$ over $\left\{|x-y|\leq dt,t\leq \tau \leq 2t\right\}$. If we take the sum over $j$ of the integrals, we obtain
\begin{align*}
   \left\|\sup_{|\cdot-y|\leq d}\|(\Psi \hat{f})\check{\ }(\cdot)|E\| \ \big|L_p\right\|+&\left\|\left(\int_0^1 t^{-sq}\sup \|(\psi(t\cdot) \hat{f})\check{\ }|E\|^q \right)^{\frac{1}{q}}\big|L_p\right\| \\
\leq c' \|(\Phi^*f)_a|L_p\|+&c'\left\|\left(\sum_{j=1}^{\infty} 2^{jsq}(\varphi_j^*f)_a)^q \right)^{\frac{1}{q}}\big|L_p\right\|
\end{align*}
and so the norms $\eqref{RychkovFolg1}$ and $\eqref{RychkovFolg2}$ are estimated from above by $c'\|f|F_{p,q}^s(E)\|$.

In the second part of the proof we want to estimate $\|f|F_{p,q}^s(E)\|$ from above again. For this we go back to step 1 of the proof of theorem $\ref{Rychkov}$, but this time interchanging the roles of $\Phi$ and $\Psi$ and of $\varphi$ and $\psi(\tau \cdot)$ in comparison to the just shown (see step 3 of the proof of theorem $\ref{Rychkov}$ for details). For given $\tau \in [1,2]$, $\Psi$ and $\psi(\tau \cdot)$ we choose functions $\Lambda^{\tau}$ and $\lambda^{\tau}=\lambda(\tau \cdot)$ (which is possible) by lemma $\ref{ExLambda}$ with the properties $\eqref{Rychkovsupport}$ and $\eqref{RychkovLambda3}$. By looking at the construction in lemma $\ref{ExLambda}$ one can see that for all $N,M \in \mathbb{N}$ there exists a $C_{N,M}$ such that
\begin{align*}
 \sup_{y \in \mathbb{R}^n} (1+|y|)^N \sum_{|\alpha|\leq M} |D^{\alpha} \Lambda^{\tau}(y)| \leq C_{N,M},
\end{align*}
i.e. $\!$that the ${\cal S}(\mathbb{R}^n)$-seminorms from $\eqref{GrundSNormen}$ can be estimated uniformly in $\tau$. This holds for $\lambda^{\tau}=\lambda(\tau \cdot)$ as well.
So we obtain the analogue of $\eqref{RychkovPsi}$, with exchanged roles,
\begin{align*}
   \left(\check{\varphi}_j*f\right)(y)=\left(\left(\check{\varphi}_j*\check{\Lambda^{\tau}}\right) * \left(\check{\Psi} * f\right)\right)(y)+ \sum_{k=1}^{\infty}\left(\left(\check{\varphi}_j*\check{\lambda^{\tau}}_k\right) * \left(\psi(\tau \cdot\right)\check{\ }_k * f)\right)(y)
\end{align*}
for all $y \in \mathbb{R}^n$. If we now follow the proof of step 1, we have to estimate the integrals from $\eqref{RychkovFaltung}$ as in $\eqref{RychkovI_j,k}$ by a constant independent of $\tau$. These are of the form
\begin{align*}
 I_{j,k}^{\tau}:=\int_{\mathbb{R}^n} \left|\left(\check{\varphi}_j*\check{\lambda^{\tau}}_k\right)(z)\right|(1+|2^kz|)^a \ dz
\end{align*}
resp.$ $ an analogue for $\Lambda^{\tau}$. To estimate the integrals from above we used lemma $\ref{RychkovI}$. Note that the constants appearing in this lemma only depend on the ${\cal S}(\mathbb{R}^n)$-seminorms and the behaviour at $0$ of $\varphi$ and $\lambda$ resp. $\!$only of the ${\cal S}(\mathbb{R}^n)$-seminorms of $\Lambda$. Hence there exists a constant independent of $\tau$ such that 
\begin{align*}
I_{j,k}^{\ \tau}\leq C_{\lambda,\varphi}\left\{
\begin{array}{l l}
 2^{(k-j)(S+1)}&, k \leq j \\
 2^{(j-k)(a-s+1)}&, k\geq j  
\end{array}
\right. .
\end{align*}
There is an analogous result for $\Lambda$. If we go on with step 1 of the proof of theorem $\ref{Rychkov}$, we get the corresponding results of $\eqref{RychkovWinner1}$ and $\eqref{RychkovWinner2}$, with exchanged roles. Hence it holds
\begin{align*}
  \sup_{y\in \mathbb{R}^n} \frac{\|\left(\check{\varphi}_j*\check{\lambda^{\tau}}_k * \psi(\tau \cdot)\check{\ }_k * f \right)(y)|E\|}{(1+2^j|x-y|)^a} &\leq  C_{\lambda,\varphi}(\psi(\tau\cdot)_k^*f)_a(x) \left\{\begin{array}{l l}
 2^{(k-j)(S+1)}&, k \leq j \\
 2^{(j-k)(-s+1)}&, k\geq j  
\end{array}
\right.
\end{align*}
and an analogue for $\Lambda^{\tau}$ with $(\Psi^*f)_a$ on the right-hand side, with $C$ independent of $\tau \in [1,2]$. Note that 
\begin{align*}
 (\psi (\tau \cdot)_{k}^*f)_a(x)&=\sup_{y\in \mathbb{R}^n} \frac{\|(\psi(2^{-k}\tau\cdot)\hat{f})\check{\ }(x-y)|E\|}{(1+2^k|y|)^a} \\
&\leq c \sup_{y\in \mathbb{R}^n} \frac{\|(\psi(2^{-k}\tau\cdot)\hat{f})\check{\ }(x-y)|E\|}{(1+2^k\tau^{-1}|y|)^a} =:(\psi_{2^{-k}\tau}^{*'}f)_a(x).
\end{align*}
With the same steps as in the proof of theorem $\ref{Rychkov}$ we obtain as a result
\begin{align*}
  2^{js}(\varphi_{j}^*f)_{a}(x) \leq C 2^{-j\delta} (\Psi^*f)_{a}(x)+C\sum_{k=1}^{\infty}  2^{ks} (\psi_{2^{-k}\tau}^{*'}f)_a(x)2^{-|j-k|\delta}
\end{align*}
and an obvious counterpart for $(\Phi^*f)_{a}$ with a certain $\delta>0$ and with $C$ independent of $\tau$. This yields
\begin{align*}
  2^{js}(\varphi_{j}^*f)_{a}(x) \leq C 2^{-j\delta} (\Psi^*f)_{a}(x)+C \left(\int_1^2 \left(\sum_{k=1}^{\infty} 2^{ks}2^{-|j-k|\delta}(\psi_{2^{-k}t}^{*'}f)_a(x)\right)^q \ \ dt\right)^{\frac{1}{q}}
\end{align*}
and an analogous result for $(\Phi^*f)_{a}(x)$. By this and a typical Minkowski/Hölder argument we obtain
\begin{align*}
 2^{js}(\varphi_{j}^*f)_{a}(x) \leq C' 2^{-j\delta_0} (\Psi^*f)_{a}(x)+C' \sum_{k=1}^{\infty} 2^{-|j-k|\delta_0}2^{ks} \left(\int_1^2 \left( (\psi_{2^{-k}t}^{*'}f)_a(x)\right)^q \ \ dt\right)^{\frac{1}{q}}
\end{align*}
for a suitable $\delta_0>0$ and an analogous result for $(\Phi^*f)_{a}(x)$ again. Now we use lemma $\ref{Delta}$ as in the proof of theorem $\ref{Rychkov}$, but this time with
\begin{align*}
 g_k:=2^{ks} \left(\int_1^2 \left( (\psi_{2^{-k}t}^{*'}f)_a(x)\right)^q \ \ dt\right)^{\frac{1}{q}}
\end{align*}
 and conclude
\begin{align*}
  \|(\Phi^*f)_{a}&|L_p\|+\left\|\left(\sum_{j=1}^{\infty} 2^{jsq}\left((\varphi_j^*f)_a\right)^q \right)^{\frac{1}{q}}\big|L_p\right\|\\
&\leq C \|(\Psi^*f)_{a}|L_p\|+C\left\|\left(\sum_{j=1}^{\infty} 2^{jsq}\int_1^2 \left( (\psi_{2^{-j}t}^{*'}f)_a\right)^q dt \right)^{\frac{1}{q}}\big|L_p\right\| \\
&\leq C'\|(\Psi^*f)_{a}|L_p\|+C'\left\|\left(\int_0^1 t^{-sq} \left((\psi_t^{*'}f)_a\right)^q \frac{dt}{t}\right)^{\frac{1}{q}}\big|L_p\right\|. 
\end{align*} 
Now we modify step 2 by applying it to $\psi(\tau \cdot)$ instead of $\psi$ for $1\leq \tau \leq 2$ to replace $(\psi_t^{*'}f)_a$ by $\psi(t\cdot)$. 

After choosing $\lambda, \Lambda \in {\cal S}(\mathbb{R}^n)$ with the desired properties \eqref{RychkovLambda4} it follows as in $\eqref{RychkovDarst}$
\begin{align}
\label{RychkovDarst2}
\begin{split}
\left(\psi(\tau\cdot)\check{\ }_j*f\right)(y)&=\left(\left[\Lambda(\tau\cdot)\check{\ }_j * \Psi(\tau\cdot)\check{\ }_j\right] * \left[\psi(\tau\cdot)\check{\ }_j*f\right]\right)(y)\\
&\hspace{3em}+\sum_{k=j+1}^{\infty}\left(\left[\psi(\tau\cdot)\check{\ }_j*\lambda(\tau\cdot)\check{\ }_k)\right] * \left[\psi(\tau\cdot)\check{\ }_k * f\right]\right)(y)
\end{split}
\end{align}
for all $y\in \mathbb{R}^n$. As there one can obtain by lemma $\ref{RychkovI}$ and in view of $k\geq j$ 
\begin{align*}
 \left|\left(\psi(\tau\cdot)\check{\ }_j*\lambda(\tau\cdot)\check{\ }_k\right)(z)\right|&=\left|2^{jn}\tau^{-n}\left(\check{\psi} * \check{\lambda}_{k-j}\right) (2^{j}\tau^{-1} z)\right| \\
&\leq C_{N,\psi,\lambda}\frac{2^{jn}\tau^{-n}2^{(j-k)N}}{(1+2^{j}\tau^{-1}|z|)^a} \leq C_{N,\psi,\lambda}'\frac{2^{jn}2^{(j-k)N}}{(1+2^{j}|z|)^a}
\end{align*}
for all $N \in \mathbb{N}$ and
\begin{align*} 
\left|\left(\Psi(\tau\cdot)\check{\ }_j*\Lambda(\tau\cdot)\check{\ }_j(z)\right)\right|&=2^{jn}\tau^{-n} \left(\check{\Psi} * \check{\Lambda}\right)(2^j\tau^{-1}z) \leq C_{\Psi,\Lambda} \frac{2^{jn}}{(1+2^j|z|)^a},
\end{align*}
where $C_{N,\psi,\lambda}'$ and $C_{\Psi,\Lambda}$ do not depend on $\tau \in [1,2]$. From $\eqref{RychkovDarst2}$ we get the analogous result of $\eqref{RychkovR}$, namely
\begin{align}
\label{RychkovFolgr1}
&(\psi_{2^{-j}\tau}^{*'}f)_a(x) \leq \! C_{N} \!\sum_{k=j}^{\infty} 2^{(j-k)N'}[(\psi_{2^{-k}\tau}^{*'}f)_a(x)]^{1-r} \!\int_{\mathbb{R}^n} \!\! \frac{2^{kn}\|\left(\psi(2^{-k}\tau\cdot)\check{\ }*f\right)(z)|E\|^r}{(1+2^k|x-z|)^{ar}} \ dz.
\end{align}
Furthermore, we derive an estimate for $\Psi$. However, this will be of slightly different shape: We start with the analogue of $\eqref{RychkovDarst2}$ for $j=1$ and $\Psi$ instead of $\psi(\tau\cdot)\check{\ }_j$
\begin{align*}
 \left(\check{\Psi}*f\right)\!(y)&\!=\!\!\left(\left[\Lambda(\tau\cdot)\check{\ } * \Psi(\tau\cdot)\check{\ }\right] * \left[\check{\Psi}*f\right]\right)\!(y)+ \sum_{k=1}^{\infty}\left(\left[\check{\Psi}*\lambda(\tau\cdot)\check{\ }_k)\right] * \left[\psi(\tau\cdot)\check{\ }_k * f\right]\right)\!(y)
\end{align*}
for all $y \in \mathbb{R}^n$. Now by lemma $\ref{RychkovI}$ we have 
\begin{align*}
  \left|\left(\check{\Psi}*\lambda(\tau\cdot)\check{\ }_k\right)(z)\right| \leq C_{\Psi,\lambda,N} \frac{\tau^{N}2^{-kN}}{(1+|z|)^a} 
\leq C_{\Psi,\lambda,N}'\frac{2^{-kN}}{(1+|z|)^a} 
\end{align*}
and obviously
\begin{align*}
  \left|\left(\Lambda(\tau\cdot)\check{\ } * \Psi(\tau\cdot)\check{\ }\right)(z)\right|\leq C_{\Psi,\Lambda}\frac{1}{(1+|z|)^a}
\end{align*}
with $C_{\Psi,\lambda,N}'$ and $C_{\Psi,\Lambda}$ independent of $\tau \in [1,2]$. Thereby we come with the same arguments as in the original proof to
\begin{align}
\label{RychkovFolgr2}
\begin{split}
 (\Psi^*f)_a&(x) \leq C_{N} [(\Psi^*f)_a(x)]^{1-r} \int_{\mathbb{R}^n} \frac{\|\left(\check{\Psi}*f\right)(z)|E\|^r}{(1+|x-z|)^{ar}} \ dz \\
&+C_{N}\sum_{k=1}^{\infty} 2^{-kN'}[(\psi_{2^{-k}\tau}^{*'}f)_{a}(x)]^{1-r} \int_{\mathbb{R}^n} \frac{2^{kn}\|\left(\psi(2^{-k}\tau\cdot)\check{\ }*f\right)(z)|E\|^r}{(1+2^k|x-z|)^{ar}}\ dz.
\end{split}
\end{align}
If $r \leq 1$, we use lemma $\ref{RychkovO}$, applied to inequality $\eqref{RychkovFolgr1}$ and $\eqref{RychkovFolgr2}$, which are valid for all $N' \in \mathbb{N}$, with $d_j=(\psi_{2^{-j}\tau}^{*'}f)_{a}(x)$. Notice that the results about condition $\eqref{RychkovLemmaO}$ can be transfered from step 3 of the proof of theorem $\ref{Rychkov}$ and hold for the $d_j$'s here, too.
Therefore, we obtain the analogue of $\eqref{RychkovAbsch}$ in the case $r\leq 1$
\begin{align}
\label{RychkovFolgRes1}
(\psi_{2^{-j}\tau}^{*'}f)_a(x)^r \leq C_{N}' \sum_{k=j}^{\infty} 2^{(j-k)Nr} \int_{\mathbb{R}^n}  \frac{2^{kn}\|\left(\psi(2^{-k}\tau\cdot)\check{\ }*f\right)(z)|E\|^r}{(1+2^k|x-z|)^{ar}} \ dz
\end{align}
and the analogue of $\eqref{RychkovAbsch4}$ for $(\Psi^*f)_a(x)^r$, where the constant $C_{N}'$ does not depend on $r \in (0,1]$, $f \in {\cal S}'(\mathbb{R}^n,E)$, $j \in \mathbb{N}$ and $\tau \in [1,2]$. As in the initial proof the assertion follows for $r>1$ as well.

In the $F_{p,q}^s(E)$-case we argue as follows: We raise to the power of $\frac{q}{r}$, integrate over $\tau \in [1,2]$ with respect to $\frac{d\tau}{\tau}$, take the $\frac{q}{r}$-th root and obtain
\begin{align*} 
&\left(\int_{2^{-j}}^{2^{-j+1}} \left((\psi_{t}^{*'}f)_a(x)\right)^q \ \frac{dt}{t}\right)^{\frac{r}{q}} \\
&\hspace{2em}\leq C_N'\left( \int_{1}^{2} \left(\sum_{k=j}^{\infty} 2^{(j-k)Nr} \int_{\mathbb{R}^n} \ \frac{2^{kn}\|\left(\psi(2^{-k}\tau\cdot)\check{\ }*f\right)(z)|E\|^r}{(1+2^k|x-z|)^{ar}} \ dz \right)^{\frac{q}{r}} \ \frac{dt}{t}\right)^{\frac{r}{q}}
\end{align*}
and 
\begin{align*}
 (\Psi^*f)_a(x)^r &\leq C_N'' \int_{\mathbb{R}^n} \frac{\|\left(\check{\Psi}*f\right)(z)|E\|^r}{(1+|x-z|)^{ar}} \ dz\\
&+C_N'' \left( \int_{1}^{2} \left(\sum_{k=1}^{\infty} 2^{-kNr} \int_{\mathbb{R}^n} \frac{2^{kn}\|\left(\psi(2^{-k}\tau\cdot)\check{\ }*f\right)(z)|E\|^r}{(1+2^k|x-z|)^{ar}}\ dz\right)^{\frac{q}{r}} \ \frac{dt}{t}\right)^{\frac{r}{q}}.
\end{align*}
If $r\leq q$, we can use the (generalized) Minkowski inequality two times and get
\begin{align*} 
&\left(\int_{2^{-j}}^{2^{-j+1}} \left((\psi_{t}^{*'}f)_a(x)\right)^q \ \frac{dt}{t}\right)^{\frac{r}{q}} \\
&\hspace{1em}\leq C_N' \sum_{k=j}^{\infty} 2^{(j-k)Nr} \int_{\mathbb{R}^n} \ \frac{2^{kn}}{(1+2^k|x-z|)^{ar}}\left( \int_{2^{-k}}^{2^{-k+1}}\|\left(\psi(t\cdot)\check{\ }*f\right)(z)|E\|^q \, \frac{dt}{t}\right)^{\frac{r}{q}} \ dz
\end{align*}
and a corresponding result for $(\Psi^*f)_a(x)^r$. This yields the estimates $\eqref{RychkovAbsch}$ and $\eqref{RychkovAbsch4}$, only with the terms
\begin{align*}
  \left(\int_{2^{-k}}^{2^{-k+1}}\|\left(\psi(t\cdot)\check{\ }*f\right)(z)|E\|^q \frac{dt}{t}\right)^{\frac{1}{q}}
\end{align*}
instead of $\|\left(\check{\psi}_k*f\right)(z)|E\|$ and
\begin{align*}
  \left(\int_{2^{-j}}^{2^{-j+1}}\left((\psi_{t}^{*'}f)_a(x)\right)^q \frac{dt}{t}\right)^{\frac{1}{q}}
\end{align*}
instead of $(\psi_j^{*}f)_a(x)$. We pick an $a$ so large that we can choose $r$ with $\frac{n}{a}<r<\min(p,q)$. Now we reconstruct the further steps in the initial proof with the given ``replacements`` and obtain immediately
\begin{align*}
 \|(\Psi^*f)_a|L_p\|&+\left\|\left(\int_0^1 t^{-sq}\left((\psi_{t}^{*'}f)_a\right)^q \frac{dt}{t} \right)^{\frac{1}{q}}\big|L_p\right\|\\ 
&\leq C'' \| (\Psi \hat{f})\check{\ }|L_p(E)\|+C''\left\|\left(\int_0^1 t^{-sq}\|(\psi(t\cdot)\hat{f})\check{\ }|E\|^q \frac{dt}{t}\right)^{\frac{1}{q}}\big|L_p\right\|.
\end{align*}
In the $B_{p,q}^s(E)$-case we start with \eqref{RychkovFolgRes1} and arrive at 
\begin{align*}
2^{jsr}(\psi_{2^{-j}\tau}^{*'}f)_a(x)^r \leq C \sum_{k=j}^{\infty} 2^{(j-k)\delta} 2^{ksr}M\left(\|\left(\psi(2^{-k}\tau\cdot)\check{\ }*f\right)|E\|^r\right)(x)
\end{align*}
as in the original proof and at an analogous result for $(\Psi^*f)_a(x)^r$. Now we take the $L_{\frac{p}{r}}$-norm, use the Minkowski inequality and the boundedness of the maximal operator from $L_{\frac{p}{r}}$ to $L_{\frac{p}{r}}$ and come to
\begin{align*}
2^{jsr}\|(\psi_{2^{-j}\tau}^{*'}f)_a(x)|L_p\|^r \leq C' \sum_{k=j}^{\infty} 2^{(j-k)\delta} 2^{ksr}\left\|\left(\psi(2^{-k}\tau\cdot)\check{\ }*f\right)|L_p(E)\right\|^r
\end{align*}
and to an analogous result for $(\Psi^*f)_a(x)^r$. Now we integrate over $\tau \in [1,2]$ with respect to $\frac{d\tau}{\tau}$, argue as in the $F_{p,q}^s(E)$-case and use a suitable estimate for the $l_q$-norm as in lemma $\ref{Delta}$. Then we obtain the desired result for $B_{p,q}^s(E)$. 

The characterizations $\eqref{RychkovFolgCharak1}$ and $\eqref{RychkovFolgCharak2}$ hold true by the same arguments as in step 4 of the proof of theorem $\ref{Rychkov}$. One just has to notice that lemma $\ref{RychkovO}$ is applied to $d_j=(\psi_{2^{-j}\tau}^{*'}f)_{a}(x)$ instead of $(\psi_j^{*}f)_{a}(x)$ which makes no big difference for the calculations.
\end{bew}

\subsection{Explicit norms and characterizations}
Below we will take a look at some examples of equivalent norms and characterizations following from theorem $\ref{Rychkov}$ which will be of use later on.	
\begin{Beispiel}                                                           
\label{RychkovLokal}  
Let $k_0,k^0 \in {\cal S}(\mathbb{R}^n)$, $\hat{k_0}(0)\neq 0$, $\hat{k^0}(0)\neq 0$. Let $N \in \mathbb{N}_0$ with $2N>s$. Then the functions $\Psi:=\hat{k_0}$ and $\psi:=\widehat{\Delta^N k^0}$ fulfil the conditions $\eqref{Rychkov2}$ and $\eqref{Rychkov3}$ for a suitable small $\varepsilon>0$ and $\eqref{Rychkov4}$ with $S=2N-1$ since $\widehat{\Delta k^0}(x)=-|x|^2 \hat{k^0}(x)$. In particular, $k_0$ and $k^0$ can be chosen such that $supp \ k_0$, $supp \ k^0 \subset B$. This is where the expression ``local means`` comes from because if $f \in {\cal S}'(\mathbb{R}^n,E)$ is e.g. a regular distribution, then we have
\begin{align*}
  (\hat{k_0}(2^{-j} \cdot) \hat{f})\check{\ }(x)= 2^{jn} \left(k_0(2^j \cdot) * f\right)(x)=\int_{B} k_0(y) f(x-2^{-j}y) \ dy
\end{align*}
and an analogue for $\Delta^N k^0$ so that for a calculation only the values of $f$ in the small domain $\left\{y \in \mathbb{R}^n:|y|\leq2^{-j}\right\}$ are necessary to know. If we set $k^N:=\Delta^N k^0$, it follows
\begin{Satz}
\label{RychkovLokalSatz}
Let $2N>s$.

(i) Let $0<p\leq\infty$ and $0<q\leq \infty$. Then
\begin{align*}
  \|f|B_{p,q}^s(E)\|_{k_0,k^N}:= \|k_0*f|L_p(E)\|+\left(\sum_{j=1}^{\infty} 2^{jsq}\|k_j^N*f|L_p(E)\|^q\right)^{\frac{1}{q}}
\end{align*}
(modified for $q=\infty$) is an equivalent norm for $\|\cdot|B_{p,q}^s(E)\|$. It holds
\begin{align*}
B_{p,q}^s(E)=\left\{f \in {\cal S}'(\mathbb{R}^n,E): 
 \|f|B_{p,q}^s(E)\|_{k_0,k^N}<\infty \right\}.
\end{align*}
(ii) Let $0<p<\infty$ and $0<q\leq \infty$. Then
\begin{align*}
  \|f|F_{p,q}^s(E)\|_{k_0,k^N}:=\|k_0*f|L_p(E)\|+\left\|\left(\sum_{j=1}^{\infty} 2^{jsq} \|k_j^N*f|E\|^q \right)^{\frac{1}{q}}\big|L_p\right\|
\end{align*}
(modified for $q=\infty$) is an equivalent norm for $\|\cdot|F_{p,q}^s(E)\|$. It holds
\begin{align*}
 F_{p,q}^s(E)=\left\{f \in {\cal S}'(\mathbb{R}^n,E): 
 \|f|F_{p,q}^s(E)\|_{k_0,k^N}<\infty \right\}.
\end{align*}
\end{Satz}
\end{Beispiel}

\begin{Beispiel}
\label{RychkovHarmonisch}
In our remarks we follow \cite{Tri97}, section 12.2, p. 59. Let $h(x):=(1+|x|^2)^{-\frac{n+1}{2}}$. By using \cite{StW90}, theorem 1.14, p. 6, we obtain 
\begin{align}
\label{RychkovFouriertransform}
  \widehat{e^{-t|\cdot|}}(x) = d_n h_{1/t}(x)=d_n t^{-n} \frac{1}{(1+|\frac{x}{t}|^2)^{\frac{n+1}{2}}}=d_n\frac{t}{(t^2+|x|^2)^{\frac{n+1}{2}}}
\end{align}
for a suitable constant $d_n>0$. The function $P(x,t)=d_n(t^2+|x|^2)^{\frac{n+1}{2}}$ and also its partial derivatives with respect to $t$ are harmonic in the domain $\left\{(x,t)\!:\! x \in \mathbb{R}^n, t>\!0\right\}$. 

Let $f \in {\cal S}(\mathbb{R}^n,E)$. Then 
\begin{align*}
 u(x,t):=(e^{-t|\cdot|} \hat{f})\check{\ }(x)= d_n \left(f *  \frac{t}{(|\cdot|^2+t^2)^{\frac{n+1}{2}}}\right)(x)
\end{align*}
is harmonic and hence also its partial derivatives with respect to $t$ given by 
\begin{align*}
 \frac{\partial^k u(x,t)}{\partial t^k}=((-1)^k|\cdot|^ke^{-t|\cdot|} \hat{f})\check{\ }(x).
\end{align*}

We choose a $\phi \in {\cal S}(\mathbb{R}^n)$ with $\phi(x)=1$ for $|x|\leq 1$ and $\phi(x)=0$ for $|x|>\frac{3}{2}$ and set $\Psi:=\phi$ and $\psi^k(\xi):=|\xi|^ke^{-\xi}$. It follows 
\begin{align*}
 (\psi^k(t\cdot) \hat{f})\check{\ }(x)= (-1)^k t^k\frac{\partial^k u(x,t)}{\partial t^k}
\end{align*}
at least for $f \in {\cal S}(\mathbb{R}^n,E)$. The functions $\Psi$ und $\psi^k$ fulfil the support conditions $\eqref{Rychkov2}$ and $\eqref{Rychkov3}$ but $\psi^k$ is not arbitrarily often differentiable. But, for instance, there exist all partial derivatives of $\psi^{2k}$ up to the order $2k$. The moment conditions $\eqref{Rychkov4}$ are fulfilled for all derivatives up to the order $2k-1$.

This means that we cannot apply theorem $\ref{Rychkov}$ including the subsequent proposition $\ref{RychkovFolgerung}$ directly for this functions. Nevertheless, we will try to obtain the desired result in this case as well. For this we will choose $k\geq k_0$ in a suitable dependence of $p,q$ and $s$.

Let $f \in B_{p,q}^s(E)$ resp.$ $ $F_{p,q}^s(E)$ be given. In consequence of $|\cdot|^ke^{-t|\cdot|} \notin {\cal S}(\mathbb{R}^n)$ we cannot use the initial definition of $(e^{-t|\cdot|} \hat{f})\check{\ }$. So we decompose $f$ into $f_1:=(\phi \hat{f})\check{\ }$ and $f_2:=((1-\phi)\hat{f})\check{\ }$. By the assumptions $f_1 \in L_p^{B_2}(E)$ and therefore by Nikolskii's inequality (see $\eqref{GrundNikolskij}$) $f_1 \in L_{\infty}$. Hence 
\begin{align*}
 u_1(x,t)=d_n \left(f_1 *  \frac{t}{(|\cdot|^2+t^2)^{\frac{n+1}{2}}}\right)(x)
\end{align*}
is well-defined, even bounded (also in $t$) as a convolution of an $L_1$-function with an $L_{\infty}(E)$-function. Moreover, the function is harmonic in the domain $\left\{(x,t): x \in \mathbb{R}^n, t>0\right\}$.
Then $u_1(x,t)$ is arbitrarily often differentiable with respect to $x$. Moreover, the (classical) partial derivatives with respect to $t$ exist and are harmonic because of $u_1(x,t)=d_n \left(h_{1/t}*f_1\right)(x)$.

The functions $[1+|\cdot|^2]^{\sigma}e^{-t|\cdot|}(1-\Phi)$, whose ${\cal S}(\mathbb{R}^n)$-seminorms for $t>\delta$ are uniformly bounded, are Fourier multipliers for $B_{p,q}^s(E)$ resp.$ $ $F_{p,q}^s(E)$ (see $\eqref{GrundFourierMult}$) for all $\sigma \in \mathbb{R}$ since $e^{-t|\cdot|}(1-\phi) \in {\cal S}(\mathbb{R}^n)$. So it follows from the lift property of these spaces (see $\eqref{GrundLift}$) and the Sobolev embeddings $\eqref{GrundSobolevB}$ and $\eqref{GrundSobolevF}$ that 
\begin{align*}
u_2(\cdot,t)=\left([1+|\cdot|^2]^{\sigma}e^{-t|\cdot|}(1-\phi)\left([1+|\cdot|^2]^{-\sigma}\hat{f}\right)\right) ^{\vee} \in B_{\infty,\infty}^s 
\end{align*}
for all $s \in \mathbb{R}$. So $u_2(\cdot,t)$ is arbitrarily often differentiable with respect to $x$ and bounded in the domain $\left\{x \in \mathbb{R}^n, t>\delta\right\}$ for a fixed $\delta>0$ (by $\eqref{GrundCub}$). The differentiability relative to $t$ is obvious. The function is harmonic by basic properties of the Fourier transform. 

Analogous assertions hold true for $|\cdot|^ke^{-t|\cdot|}$ instead of $e^{-t|\cdot|}$ and therefore for the partial derivatives of $u_2(x,t)$. So $u(x,t):=u_1(x,t)+u_2(x,t)$ is well-defined for all $f \in B_{p,q}^s(E)$ resp.$ $ $F_{p,q}^s(E)$ for arbitrary $s$, $0<p\leq \infty$ (resp.$ $ $<\infty$) and $0<q\leq \infty$, arbitrarily often differentiable, harmonic in the domain $\left\{ x \in \mathbb{R}^n, t>0\right\}$ and bounded on $\left\{x \in \mathbb{R}^n, t>\delta\right\}$ for fixed $\delta>0$. An analogue is valid for the partial derivatives.

Now we have the necessary tools together to formulate and proof the desired result.
\begin{Satz}
\label{RychkovHarmSatz}
Let $d>0$ and $s \in \mathbb{R}$.

(i) Let $0<p\leq \infty$, $0<q \leq \infty$. Then there exists a $k_0$ such that for all $k \geq k_0$ 
\begin{align*}
  \|f|B_{p,q}^s(E)\|_{\phi,k}:=\|(\phi \hat{f})\check{\ }|L_p(E)\|+\left(\int_0^{1} t^{(k-s)q}\|\frac{\partial^k u(\cdot,t)}{\partial t^k}|L_p(E)\|^q \ \frac{dt}{t} \right)^{\frac{1}{q}}
\end{align*}
and 
\begin{align*}
  \|f|B_{p,q}^s(E)\|_{\phi,k}^{\sup}:=&\left\|\sup_{|\cdot-y|\leq d}\|(\phi \hat{f})\check{\ }(y)|E\| \big|L_p\right\| \!+\!\left(\int_0^{1} t^{(k-s)q}\left\|\sup \| \frac{\partial^k u(y,\tau)}{\partial t^k}|E\| \big|L_p\right\|^q \ \frac{dt}{t}\right)^{\frac{1}{q}}
\end{align*}
are equivalent norms for $\|\cdot|B_{p,q}^s(E)\|$, where $\sup$ is the supremum for a fixed $x \in \mathbb{R}^n$ over $\left\{|x-y|\leq dt,t\leq \tau \leq 2t\right\}$. It holds
\begin{align*}
 B_{p,q}^s(E)=\left\{f \in \mathscr{C}^{-\infty}(E): \|f|B_{p,q}^s(E)\|_{\phi,k}<\infty \right\}
\end{align*}
 and 
\begin{align*}
 B_{p,q}^s(E)=\left\{f \in \mathscr{C}^{-\infty}(E): 
\|f|B_{p,q}^s(E)\|_{\phi,k}^{\sup} <\infty \right\}.
\end{align*}

(ii) Let $0<p< \infty$, $0<q \leq \infty$. Then there exists a $k_0$ such that for all $k \geq k_0$ 
\begin{align*}
  \|f|F_{p,q}^s(E)\|_{\phi,k}:=\left\|(\phi \hat{f})\check{\ }\big|L_p(E)\right\|+\left\|\left(\int_0^1 t^{(k-s)q} \|\frac{\partial^k u(\cdot,t)}{\partial t^k}|E\|^q \ \frac{dt}{t} \right)^{\frac{1}{q}}\big|L_p\right\|
\end{align*}
and
\begin{align*}
\|f|F_{p,q}^s(E)\|_{\phi,k}^{\sup}:=&\left\|\sup_{|\cdot-y|\leq d}\|(\phi \hat{f})\check{\ }(y)|E\| \big|L_p\right\| +\left\|\left(\int_0^1 t^{-sq}\sup \|\frac{\partial^k u(y,\tau)}{\partial t^k}|E\|^q \right)^{\frac{1}{q}}\big|L_p\right\|
\end{align*}
are equivalent norms for $\|\cdot|F_{p,q}^s(E)\|$, where $\sup$ is the supremum for a fixed $x \in \mathbb{R}^n$ over $\left\{|x-y|\leq dt,t\leq \tau \leq 2t\right\}$. It holds
\begin{align}
\label{RychkovHarmCh1}
 F_{p,q}^s(E)=\left\{f \in \mathscr{C}^{-\infty}(E): \|f|F_{p,q}^s(E)\|_{\phi,k}<\infty \right\}
\end{align}
 and 
\begin{align}
\label{RychkovHarmCh2}
 F_{p,q}^s(E)=\left\{f \in \mathscr{C}^{-\infty}(E): 
\|f|F_{p,q}^s(E)\|_{\phi,k}^{\sup} <\infty \right\}.
\end{align}
\end{Satz}
\begin{bew}
We like to recall the relation
\begin{align*}  
(-1)^kt^k\frac{\partial^k u(x,t)}{\partial t^k}=\left((t|\cdot|)^ke^{-t|\cdot|} \hat{f}\right)^{\vee}(x).
\end{align*} 
That explains the form of the norm in this theorem in our context. The proof is a step-by-step repetition of theorem $\ref{Rychkov}$ resp.$ $ proposition $\ref{RychkovFolgerung}$, where we use that the function  $\psi:=|\cdot|^k e^{-|\cdot|}$ behaves like a ${\cal}S(\mathbb{R}^n)$-function away from $0$ and fulfils as many moment conditions as we want if we only choose $k$ large enough.

We have to say one word about condition $\eqref{RychkovLemmaO}$. It holds
\begin{align*}
 d_j=(\psi_j^*f)_a(x)\leq \left(((1-\phi)\psi)_j^*f\right)_a(x)+((\phi\psi)_j^*f)_a(x).
\end{align*}
The first summand can be estimated as in step 3 of the proof of theorem $\ref{Rychkov}$ because $(1-\phi)\psi \in {\cal S}(\mathbb{R}^n)$. For the second summand we have
\begin{align*}
 \|(\phi\psi)\check{\ }_j*f|L_{\infty}(E)\|\leq \|\check{\phi}_j*f|L_{\infty}\|\cdot \|\check{\psi}_j|L_1\|\leq c \|\check{\phi}_j*f|L_{\infty}(E)\|.
\end{align*}
We choose a smooth dyadic resolution of unity consisting of the functions $\{\varphi_j\}_{j \in \mathbb{N}}$ (see $\eqref{GrundResolution}$) and obtain by Nikolskii's inequality (see $\eqref{GrundNikolskij}$)
\begin{align*}
 \|\check{\phi}_j*f|L_ {\infty}(E)\| & \leq\sum_{k=0}^{j+1} \|\check{\phi}_j*(\varphi_k\hat{f})\check{\ }|L_ {\infty}(E)\|\leq c \sum_{k=0}^{j+1} \|(\varphi_k\hat{f})\check{\ }|L_ {\infty}(E)\|\\
&\leq c' \sum_{k=0}^{j+1} 2^{\frac{kn}{p}} \|(\varphi_k\hat{f})\check{\ }|L_p(E)\|\leq c'' \left(\sum_{k=0}^{j+1} 2^{k\left(\frac{n}{p}+\varepsilon\right)q} \|(\varphi_k\hat{f})\check{\ }|L_p(E)\|^q\right)^{\frac{1}{q}} \\
&\leq c' \max \left(2^{j\left(\frac{n}{p}+\varepsilon-s\right)},1\right) ||f|B_{p,q}^s(E)||.
\end{align*}
So the desired condition $\eqref{RychkovLemmaO}$ with $N_0\geq\lceil\frac{n}{p}+\varepsilon-s\rceil$ follows.

Because we assumed $f \in \mathscr{C}^{-\infty}(E)$ for defining the convolution with $e^{-|\cdot|}$ a priori, we obtain the best possible results with the characterizations $\eqref{RychkovHarmCh1}$ and $\eqref{RychkovHarmCh2}$.
\end{bew}

\begin{Bemerkung}
We follow \cite{Tri97}, theorem 12.5 (i) and (ii), p. 64. We want to replace $\phi$ in our proposition $\ref{RychkovHarmSatz}$ by the function $e^{-|\cdot|}$ so that $(\Psi \hat{f}) \check{\ }$ is harmonic as well. But this will only work for $p>\frac{n}{n+1}$. Namely in this case $m(\xi)=e^{-|\xi|}$ is a Fourier multiplier for $L_p^{B}(E)$, i.e. that there exists a $C>0$ such that it holds
\begin{align*}
 \|(e^{-|\cdot|} \hat{f})\check{\ }|L_p(E)\| \leq C \|f|L_p(E)\|.
\end{align*}
Let's justify this: Let $0<p<1$ and $\lambda>n\left(\frac{1}{p}-1\right)$. Then it follows from $\eqref{GrundFaltung2}$ that there exists a constant $c>0$ such that for all $m:\mathbb{R}^n \rightarrow \mathbb{C}$ with $\check{m} \in L_1$ and for all $f \in L_p^{B}(E)$
\begin{align*}
 \|(m \hat{f})\check{\ }|L_p(E)\|&=\|(m\phi \hat{f})\check{\ }|L_p(E)\| \\
&\leq \|\left(m\phi\right)\check{\ }|L_p\| \cdot \|f|L_p(E)\| \\
 &\leq c\|(1+|\cdot|^2)^{\frac{\lambda}{2}} \left(m\phi\right)\check{\ }|L_1\| \cdot \|f|L_p(E)\| \\
 &\leq c\|(1+|\cdot|^2)^{\frac{\lambda}{2}} \check{m}|L_1\| \cdot \|f|L_p(E)\|.
\end{align*}
If $p\geq 1$, then we have by $\eqref{GrundFaltung1}$
\begin{align*} 
\|(m \hat{f})\check{\ }|L_p(E)\| &\leq \|\check{m}|L_1\| \cdot \|f|L_p(E)\|.
\end{align*}
The terms on the right-hand side are precisely finite if and only if $p>\frac{n}{n+1}$, see $\eqref{RychkovFouriertransform}$.

Using
\begin{align*}
\left(|\cdot|^k e^{-|\cdot|}\right)\hat{\ }=(-1)^k d_n\cdot\frac{\partial^k h_{1/t}(x,1)}{\partial t^k}
\end{align*}
and the structure of the function $h$ (see $\eqref{RychkovFouriertransform}$) it follows by the same arguments that the functions $|\cdot|^k e^{-|\cdot|}$ are Fourier multipliers for $L_p^B(E)$ with $p>\frac{n}{n+1}$ as well.
\end{Bemerkung}

\begin{Satz}
\label{RychkovHarmLast}
Let $d>0$ and $s \in \mathbb{R}$

(i) Let $\frac{n}{n+1}<p\leq \infty$, $0<q \leq \infty$ and $f \in B_{p,q}^s(E)$. Then there exist a $k_0 \in \mathbb{N}$ and a $c>0$ such that for all $k \in \mathbb{N}$ with $k \geq k_0$ 
\begin{align*}
  \|u(\cdot,1)|L_p(E)\|+\left(\int_0^{1} t^{(k-s)q}\|\frac{\partial^k u(\cdot,t)}{\partial t^k}|L_p(E)\|^q \ \frac{dt}{t} \right)^{\frac{1}{q}}
\end{align*}
and 
\begin{align*}
  \sum_{l=0}^{k-1} \left\|\sup_{|\cdot-y|\leq d}\|\frac{\partial^l u(y,1)}{\partial t^l}|E\| \ \big|L_p\right\|+\left(\int_0^{1} t^{(k-s)q}\left\|\sup \| \frac{\partial^k u(y,\tau)}{\partial t^k}|E\|\big|L_p\right\|^q \ \frac{dt}{t}\right)^{\frac{1}{q}}
\end{align*}
are bounded from above by $c\, \|\cdot|B_{p,q}^s(E)\|$ where $\sup$ is the supremum for a fixed $x \in \mathbb{R}^n$ over $\left\{|x-y|\leq dt, t\leq \tau \leq 2t\right\}$.

(ii) Let $\frac{n}{n+1}<p< \infty$, $0<q \leq \infty$ and let $f \in F_{p,q}^s(E)$. Then there exist a $k_0 \in \mathbb{N}$ and a $c>0$ such that for all $k \in \mathbb{N}$ with $k \geq k_0$ 
\begin{align*}
  \|u(\cdot,1)|L_p(E)\|+\left\|\left(\int_0^1 t^{(k-s)q} \|\frac{\partial^k u(\cdot,t)}{\partial t^k}|E\|^q \ \frac{dt}{t} \right)^{\frac{1}{q}}\big|L_p\right\|
\end{align*}
and
\begin{align}
\label{RychkovBspNorm4}
  \sum_{l=0}^{k-1} \left\|\sup_{|\cdot-y|\leq d}\|\frac{\partial^l u(y,1)}{\partial t^l}|E\| \ \big|L_p\right\|+\left\|\left(\int_0^1 t^{(k-s)q}\sup \|\frac{\partial^k u(y,\tau)}{\partial t^k}|E\|^q \right)^{\frac{1}{q}}\big|L_p\right\|
\end{align}
are bounded from above by $c\, \|\cdot|F_{p,q}^s(E)\|$ where $\sup$ is the same as in part (i).
\end{Satz}
\begin{bew}
If we take a look at proposition $\ref{RychkovHarmSatz}$ proven before, it suffices to estimate the first term of $\eqref{RychkovBspNorm4}$ by $||f|F_{p,q}^s(E)||$ (in the $F_{p,q}^s(E)$-case). 

Let $\phi$ be chosen as before and $\Psi(\xi):=\phi(\xi) |\xi|^l e^{-|\xi|}$. We obtain
\begin{align}
\label{RychkovHarmAbsch0}
\begin{split}
 \sup_{|x-y|\leq d}\|&\frac{\partial^l u(y,1)}{\partial t^l}|E\| \\
&\leq \sup_{|x-y|\leq d} \|(\phi |\cdot|^le^{-|\cdot|} \hat{f})\check{\ }(y)|E\|+\sup_{|x-y|\leq d} \|((1-\phi) |\cdot|^l e^{-|\cdot|} \hat{f})\check{\ }(y)|E\| \\
&\leq c(\Psi^*f)_a(x) + c\sup_{|x-y|\leq d} \|((1-\phi) e^{-|\cdot|} \hat{f})\check{\  }(y)|E\|
\end{split}
\end{align}
for all $a>0$, where $c$ depends on $a$. By $\eqref{GrundPeetre}$ it follows if we choose $a>\frac{n}{p}$ 
\begin{align*}
\|(\Psi^*f)_a(x)|L_p\| \leq c' \|(\Psi \hat{f})\check{\ }|L_p(E)\|. 
\end{align*}
Now we use that $|\cdot|^le^{-|\cdot|}$ are Fourier multipliers for $L_p^B(E)$ with $p>\frac{n}{n+1}$ by the remark before and obtain
\begin{align}
\label{RychkovHarmAbsch1}
\|(\Psi^*f)_a(x)|L_p\| \leq c'' \|(\phi \hat{f})\check{\ }|L_p(E)\|. 
\end{align}
Furthermore, let $g \in B_{p,1}^{\frac{n}{p}+\varepsilon}(E)$ for an $\varepsilon>0$. Then we have by \cite{Tri92}, remark 1, p.$ $ 128
\begin{align*}
 \|\sup_{|\cdot-y|\leq d} \|g|E\||L_p\| &\leq c\|g|B_{p,1}^{\frac{n}{p}+\varepsilon}(E)\|.
\end{align*}
Hereby we obtain for $g:=\left((1-\phi) |\cdot|^l e^{-|\cdot|} \hat{f}\right)\check{\ }$ with $(1-\phi) |\cdot|^l e^{-|\cdot|} \in {\cal S}(\mathbb{R}^n)$
\begin{align}
\label{RychkovHarmAbsch2}
\begin{split} 
\left\|\sup_{|x-y|\leq d} \|((1-\phi)|\cdot|^l e^{-|\cdot|} \hat{f})\check{\  }(y)|E\|\ |L_p\right\| &\leq \|((1-\phi) |\cdot|^le^{-|\cdot|} \hat{f})\check{\  }|B_{p,1}^{1+\frac{n}{p}}(E)\| \\
&\leq c \|f|F_{p,q}^s(E)\|.
\end{split}
\end{align}
If we put both results
$\eqref{RychkovHarmAbsch1}$ and $\eqref{RychkovHarmAbsch2}$ into $\eqref{RychkovHarmAbsch0}$, we arrive at the desired estimate
\begin{align*}
\left\|\sup_{|x-y|\leq d}\|\frac{\partial^l u(y,1)}{\partial t^l}\big|E\right\|&\leq c \|(\phi \hat{f})\check{\ }|L_p(E)\|+\|f|F_{p,q}^s(E)\| \\
 &\leq c'\|f|F_{p,q}^s(E)\|.
\end{align*} 
\end{bew}
\end{Beispiel}

	\section{Atomic characterizations of vector-valued function spaces}
	\subsection{Atomic and harmonic representations}
 After dealing with the necessary arrangements we now take a look at atomic representations. It is our aim to represent every element of a function space $B_{p,q}^s(E)$ resp. $F_{p,q}^s(E)$ as a preferably easy (infinite) linear combination of ``good-natured`` functions. To this we describe the concept of atoms as one can find it in \cite{Tri97}, definition 13.3, p. 73. Thereby $Q_{\nu,m}:=\{x \in \mathbb{R}^n: |x_i-2^{-\nu}m_i|\leq 2^{-\nu-1}\}$ stands for the cube with sides parallel to the axes and with the center at $2^{-\nu}m$ and side length $2^{-\nu}$ for $m \in \mathbb{Z}^n$ and $\nu \in \mathbb{N}_0$.
\begin{Definition}
\label{Atoms}
(i) Let $K \in \mathbb{N}_0$ and $d>1$. A $K$ times differentiable (in the case $K=0$ continuous) function $a:\mathbb{R}^n \rightarrow E$ is called ($E$-valued) $1$-atom (more exactly $1_K$-atom) if 
\begin{align*}
 supp \ a &\subset d \cdot Q_{0,m} \text{ for an } m \in \mathbb{Z}, \\
 \|D^{\alpha} a(x)|E\| &\leq 1 \text{ for all } |\alpha|\leq K.
\end{align*}

 (ii) Let $s \in \mathbb{R}$, $0<p\leq\infty$, $K \in \mathbb{N}_0$, $L+1 \in \mathbb{N}_0$ and $d>1$.
A $K$ times differentiable (in the case $K=0$ continuous) function $a:\mathbb{R}^n \rightarrow E$ is called ($E$-valued) $(s,p)$-atom (more exactly $(s,p)_{K,L}$-atom) if there exists a $\nu \in \mathbb{N}_0$ such that
\begin{align}
 \label{Atom1}supp \ a &\subset d \cdot Q_{\nu,m} \text{ for an } m \in \mathbb{Z},  \\
 \label{Atom2}\|D^{\alpha} a(x)|E\| &\leq 2^{-\nu\left(s-\frac{n}{p}\right)+|\alpha|\nu} \text{ for all } |\alpha|\leq K, \\
 \label{Atom3}\int_{\mathbb{R}^n} &x^{\beta} a(x) \ dx=0 \text{ for all }  |\beta|\leq L.
\end{align}
\end{Definition}
In particular,  $a_{\nu,m}e_{\nu,m}$ is a vector-valued $(s,p)_{K,L}$-atom if $a_{\nu,m}$ is a scalar (i.e. $\mathbb{C}$-valued) $(s,p)_{K,L}$-atom and $e_{\nu,m} \in U_E=\left\{x \in E: \|x|E\|=1\right\}$ .

Furthermore, we introduce the sequence spaces $b_{p,q}$ and $f_{p,q}$ whose use will become clear in the following. At this we refer to \cite{Tri97}, definition 13.5, p. 74.
\begin{Definition}
 Let $0<p\leq \infty$, $0<q\leq \infty$ and 
 \begin{align*}
  \lambda=\left\{ \lambda_{\nu,m} \in \mathbb{C}: \nu \in \mathbb{N}_0, m \in \mathbb{Z}^n\right\}.
 \end{align*}
In addition, let
\begin{align*}
 b_{p,q}:=\left\{\lambda: \|\lambda|b_{p,q}\|=\left(\sum_{\nu=0}^{\infty} \left(\sum_{m \in \mathbb{Z}^n} |\lambda_{\nu,m}|^p\right)^{\frac{q}{p}} \right)^{\frac{1}{q}} <\infty \right\}
\end{align*}
and
\begin{align*}
 f_{p,q}:=\left\{\lambda: \|\lambda|f_{p,q}\|=\left\|\left(\sum_{\nu=0}^{\infty} \sum_{m \in \mathbb{Z}^n} |\lambda_{\nu,m}\chi_{\nu,m}^{(p)}(\cdot)|^q\right)^{\frac{1}{q}}\big|L_p\right\| <\infty \right\}
\end{align*}
(modified in the cases $p=\infty$ or $q=\infty$), where $\chi_{\nu,m}^{(p)}$ is the $L_p$-normalized characteristic function of the cube $Q_{\nu,m}$.
\end{Definition}
The following lemma is oriented towards \cite{Tri97}, corollary 13.9, p. 81 which considers the scalar case. But we modify the original proof a bit. Here we get a first clue how the sought representation of all functions from $B_{p,q}^s(E)$ resp. $F_{p,q}^s(E)$ looks like.
\begin{Lemma}
\label{HarmS'-KonvAtom}
Let $0<p\leq\infty$ resp. $<\infty$, $0<q\leq \infty$ and $s \in \mathbb{R}$. Let $K \in \mathbb{N}_0$, $L+1 \in \mathbb{N}_0$ with
\begin{align}
  \label{KonvSBed1} K\geq 1+\lfloor s\rfloor \text{ and } \ L\geq \lfloor\sigma_p-s\rfloor.
\end{align}
Then 
\begin{align*}
 \sum_{\nu=0}^{\infty} \sum_{m \in \mathbb{Z}^n} \lambda_{\nu,m} a_{\nu,m}(x)
\end{align*}
converges unconditionally in ${\cal S}'(\mathbb{R}^n,E)$, where $a_{\nu,m}$ are $E$-valued $1_K$-atoms (for $\nu=0$) or $E$-valued $(s,p)_{K,L}$-atoms (for $\nu \in \mathbb{N}$) and $\lambda \in b_{p,q}$ or $\lambda \in f_{p,q}$.
\end{Lemma}
\begin{bew}
Let $\varphi \in {\cal S}(\mathbb{R}^n)$. In view of $\eqref{Atom3}$ we obtain
\begin{align*}
 \int_{\mathbb{R}^n} &  \lambda_{\nu,m} a_{\nu,m}(x)  \varphi(x) \ dx =\int_{\mathbb{R}^n} \lambda_{\nu,m} a_{\nu,m}(x) \left(\varphi(x) - \sum_{|\beta|\leq L} c_{\beta}^{\nu,m}(x-2^{-\nu}m)^{\beta} \right),
\end{align*}
where $c_{\beta}^{\nu,m} \in \mathbb{C}$ is the coefficient for $\beta$ in the Taylor expansion of $\varphi$ at $2^{-\nu}m$. 
The modulus of the difference under the integral can be estimated from above (with arbitrary $M>0$) by 
\begin{align*}
  c\, 2^{-\nu(L+1)}(1+|x|^2)^{-\frac{M}{2}}& \sup_{y \in \mathbb{R}^n} (1+|y|^2)^{\frac{M}{2}} \sum_{|\gamma|\leq L+1} \! |D^{\gamma} \varphi(y)|=c \ 2^{-\nu(L+1)}(1+|x|^2)^{-\frac{M}{2}} \|\varphi\|_{M,L+1}.
\end{align*}
In the case $1\leq p \leq \infty$ we have $L+1>-s$ (by $\eqref{KonvSBed1}$) and so by using $\eqref{Atom2}$ 
\begin{align*}
 2^{-\nu(L+1)} \|a_{\nu,m}|E\| \leq 2^{-\nu \left(s-\frac{n}{p}\right)}2^{-\nu(L+1)} \leq \chi_{\nu,m}^{(p)} 2^{-\nu \varkappa}
\end{align*}
with a $\varkappa>0$. Keeping in mind that for fixed $\nu$ the supports of $a_{\nu,m}$ are ``nearly`` disjoint we obtain together with Hölder's inequality
\begin{align*}
 &\sum_{m}\left\| \int_{\mathbb{R}^n}  \lambda_{\nu,m} a_{\nu,m}(x)  \varphi(x) \ dx\big|E\right\| \\ 
&\hspace{5em}\leq c \|\varphi\|_{M,L+1} \sum_{m} \int_{\mathbb{R}^n}  2^{-\nu(L+1)} \|\lambda_{\nu,m}a_{\nu,m}(x)|E\| (1+|x|^2)^{-\frac{M}{2}} \ dx  \\
&\hspace{5em}\leq c' \|\varphi\|_{M,L+1} \left(\int_{\mathbb{R}^n} \sum_{m} \left(2^{-\nu(L+1)} \|\lambda_{\nu,m}a_{\nu,m}(x)|E\|\right)^p\right)^{\frac{1}{p}}  \\
&\hspace{5em}\leq c''2^{-\nu \varkappa} \|\varphi\|_{M,L+1} \left(\sum_{m} |\lambda_{\nu,m}|^p\right)^{\frac{1}{p}}.
\end{align*}
if we only choose $M$ so large that $Mp'>n$. Because of $\varkappa>0$ it follows  
\begin{align}
\label{KonvSForm1}
 \begin{split}
\sum_{\nu} \sum_{m}\left\|\int_{\mathbb{R}^n}   \lambda_{\nu,m} a_{\nu,m}(x)  \varphi(x) \ dx\big|E\right\|
\leq & \, c''\,\|\varphi\|_{M,L+1} \sum_{\nu} 2^{-\nu \varkappa} \left(\sum_{m}  |\lambda_{\nu,m}|^p\right)^{\frac{1}{p}} \\
 \leq & \, C \, \|\varphi\|_{M,L+1} \cdot \|\lambda_{\nu,m}|b_{p,\infty}\|.
\end{split}
\end{align}
Because of $b_{p,q} \hookrightarrow b_{p,\infty}$ we have shown the absolute convergence of the above series in the Banach space $E$. Hence the series itself converges unconditionally in the Banach space $E$. This shows the desired claim by an admissible commutation of the integral and the sums.

In the case $0<p<1$ we have $L+1>-s+\frac{n}{p}-n$ by the assumptions instead and hence
\begin{align*}
 2^{-\nu(L+1)} \|a_{\nu,m}|E\| \leq c\, 2^{-\nu \left(s-\frac{n}{p}\right)}2^{-\nu(L+1)} \leq c\, 2^{\nu n} 2^{-\nu \varkappa}.
\end{align*}
So we obtain the above estimates $\eqref{KonvSForm1}$ for $p=1$ and because of $b_{p,q} \hookrightarrow b_{1,q} \hookrightarrow b_{1,\infty}$ the convergence in ${\cal S}'(\mathbb{R}^n,E)$ for $p<1$ follows, too.

\noindent If $\lambda \in f_{p,q}$, then $\lambda \in b_{p,\infty}$ and hence the convergence in ${\cal S}'(\mathbb{R}^n,E)$ follows as well.

Furthermore note that the condition $K=0$ would have sufficed for the whole proof, i.e. taking continuous atoms with suitable boundary conditions without any restrictions on the derivatives but with possible moment conditions.
\end{bew}

The next proposition gives a characterization of such sums as elements of the function spaces $B_{p,q}^s(E)$ and $F_{p.q}^s(E)$. At this we stick to \cite{Tri97}, theorem 13.8, p. 75, step 2, which treats the scalar case. In the last part of the proof we give a slight modification due to a small gap regarding the maximal function in the original proof. Otherwise the proof can be taken over nearly verbatim.
\begin{Satz}
\label{HarmIfAtom}
 (i) Let $0<p\leq \infty$, $0<q\leq \infty$ and $s \in \mathbb{R}$. Let $K \in \mathbb{N}_0$ and $L+1 \in \mathbb{N}_0$ with
\begin{align*}
K\geq 1+\lfloor s\rfloor \ and \ L\geq \lfloor\sigma_p-s\rfloor.
\end{align*}
Then every $f \in {\cal S}'(\mathbb{R}^n,E)$ which can be represented by
\begin{align*}
 f=\sum_{\nu \in \mathbb{N}_0} \sum_{m \in \mathbb{Z}^n} \lambda_{\nu,m} a_{\nu,m}
\end{align*}
in ${\cal S}'(\mathbb{R}^n,E)$ belongs to $B_{p,q}^s(E)$. Thereby $a_{\nu,m}$ are $E$-valued $1_K$-atoms (for $\nu=0$) or $E$-valued $(s,p)_{K,L}$-atoms (for $\nu \in \mathbb{N}$) and $\lambda \in b_{p,q}$.
Furthermore, there exists a constant $c$ independent of $f$, $\lambda$ and $a_{\nu,m}$, i.e. independent of the found representation of $f$ such that
\begin{align*}
 \|f|B_{p,q}^s(E)\| \leq c \|\lambda|b_{p,q}\|.
\end{align*}
 (ii) Let $0<p<\infty$, $0<q\leq \infty$ and $s \in \mathbb{R}$. Let $K \in \mathbb{N}_0$ and $L+1 \in \mathbb{N}_0$ with
\begin{align}
  \label{HarmAtomBed2} K\geq 1+\lfloor s \rfloor \ and \ L\geq \lfloor\sigma_{p,q}-s\rfloor.
\end{align}
Then every $f \in {\cal S}'(\mathbb{R}^n,E)$ which can be represented by
\begin{align*}
 f=\sum_{\nu \in \mathbb{N}_0} \sum_{m \in \mathbb{Z}^n} \lambda_{\nu,m} a_{\nu,m}
\end{align*}
in ${\cal S}'(\mathbb{R}^n,E)$ belongs to $F_{p,q}^s(E)$. Thereby $a_{\nu,m}$ are $E$-valued $1_K$-atoms (for $\nu=0$) or $E$-valued $(s,p)_{K,L}$-atoms (for $\nu \in \mathbb{N}$) and $\lambda \in f_{p,q}$.
Furthermore, there exists a constant $c$ independent of $f$, $\lambda$ and $a_{\nu,m}$, i.e. independent of the found representation of $f$ such that
\begin{align*}
 \|f|F_{p,q}^s(E)\| \leq c \|\lambda|f_{p,q}\|.
\end{align*}
\end{Satz}
\begin{bew}
In the proof we rely on the equivalent quasi-norm from proposition $\ref{RychkovLokalSatz}$ which results from theorem $\ref{Rychkov}$. 
We choose the functions $k_0, k^0 \in {\cal S}(\mathbb{R}^n)$ and hence also $k^N:=\Delta^N k^0$ so that they have compact support, i.e. $supp \ k_0, supp\ k^0 \subset e \cdot B$ for an $e>0$. Let $a_{\nu,m}$ with $\nu \in \mathbb{N}$ and $m \in \mathbb{Z}^n$ be an $E$-valued $(s,p)_{K,L}$-atom by definition $\ref{Atoms}$. If $\nu=0$, let $a_{\nu,m}$ be an $E$-valued $1_K$-atom. Then we can take over the argumentation from \cite{Tri97}, theorem 13.8, p. 75, step 2: If $j \geq \nu$, we get
\begin{align}
\label{HarmAbs1}
\begin{split}
 2^{js}\|\left( k^N_j * a_{\nu,m}\right)(x)|E\|\leq c\, 2^{-\varkappa(j-\nu)} \tilde{\chi}_{\nu,m}^{(p)}(x),
\end{split}
\end{align}
where $\tilde{\chi}_{\nu,m}^{(p)}(x)$ is the $L_p$-normalised characteristic function of the cube $c \cdot Q_{\nu,m}$ and $\varkappa>0$ by $\eqref{HarmAtomBed2}$. The case $\nu=0$, i.e. the case in which $a_{\nu,m}$ is a $1_K$-atom can be treated in the same way.

Let now $j< \nu$. We obtain from \cite{Tri97}, theorem 13.8, p. 75, step 2 that
\begin{align}
\label{HarmAbsch}
 2^{js} \|\left(k^N_j * a_{\nu,m}\right)(x)|E\|\leq c2^{j(s+n)} 2^{-\nu\left(s-\frac{n}{p}\right)} 2^{(j-\nu)(L+1)}\int_{|y|\leq e2^{-j}} \tilde{\chi}_{\nu,m}(x-y) \ dy.
\end{align}
Now the integral on the right-hand side is at most $d^n2^{-\nu n}$ and vanishes if we have $|x-2^{-\nu}m|>d2^{-\nu}+e2^{-j}$, i. e. if $x \notin c2^{\nu-j}\cdot Q_{\nu,m}$ for a suitable $c>0$, observing $j<\nu$. Altogether the integral is smaller or equal to $d^n 2^{-\nu n} \chi(c2^{\nu-j}Q_{\nu,m})(x)$. We have
\begin{align}
\label{HarmMaximalProb}
\begin{split}
 \sum_{m} |\lambda_{\nu,m}| \int_{|y|\leq e2^{-j}} \tilde{\chi}_{\nu,m}(x-y) \ dy &\leq \sum_{m}|\lambda_{\nu,m}|  d^n 2^{-\nu n} \chi(c2^{\nu-j}Q_{\nu,m})(x) \\
&=d^n 2^{-\nu n} \sum_{m\in D_x}|\lambda_{\nu,m}|,
\end{split}
\end{align}
where $D_x:=\left\{m \in \mathbb{Z}^n: \, x \in c2^{\nu-j}Q_{\nu,m} \right\}$. Then let
$$E_x:=\bigcup_{m \in D_x} c2^{\nu-j}Q_{\nu,m}.$$
There is a constant $c'>c$ independent of $m$ and $\nu$ such that $E_x \subset B_{c'2^{-j}}(x)$. Simultaneously it holds
\begin{align*}
 M\left(\sum_{m} |\lambda_{\nu,m}| \chi(Q_{\nu,m})\right)^w &\geq \frac{1}{|B_{c'2^{-j}}(x)|} \int_{B_{c'2^{-j}}(x)} \left(\sum_{m} |\lambda_{\nu,m}| \chi(Q_{\nu,m})(y)\right)^w \ dy \\
&\geq c'' 2^{(j-\nu)n}  \sum_{m \in D_x} |\lambda_{\nu,m}|^w
\end{align*}
because $Q_{\nu,m} \subset E_x$ for $m \in D_x$ and the $Q_{\nu,m}$ are pairwise disjoint. Together with \eqref{HarmMaximalProb} (observing $w<1$) this yields
\begin{align*}
  \sum_{m} |\lambda_{\nu,m}| \!\int_{|y|\leq e2^{-j}} \tilde{\chi}_{\nu,m}(x-y) \ dy &\leq C 2^{-\nu n} 2^{(\nu-j)\frac{n}{w}} \!\left(M\left(\sum_{m} |\lambda_{\nu,m}| \chi(Q_{\nu,m})\right)^w\right)^{\frac{1}{w}}(x).
\end{align*}
If we put this into $\eqref{HarmAbsch}$ and replace the characteristic function by the $L_p$-normalized characteristic function, we can conclude
\begin{align}
\label{HarmAbs2} 2^{js} &\left\|\left(k^N_j * \sum_{m} |\lambda_{\nu,m}| a_{\nu,m}\right)(x)\big|E\right\|\leq C' 2^{-(\nu-j)\varkappa} \left(M\left(\sum_{m} |\lambda_{\nu,m}| \chi_{\nu,m}^{(p)}\right)^w\right)^{\frac{1}{w}}(x).
\end{align}
Here we have $\varkappa>0$ if we choose $w$ close enough to $\min(1,p,q)$ resp. $\min(1,p)$. If we now use 
$\eqref{HarmAbs1}$ and $\eqref{HarmAbs2}$, we obtain
\begin{align*}
 2^{js} \left\|\left(k^N_j * \sum_{\nu,m} \lambda_{\nu,m} a_{\nu,m}(x)\right)\big|E\right\|&\leq c \sum_{\nu\leq j,m} 2^{-|j-\nu|\varkappa}|\lambda_{\nu,m}| \tilde{\chi}_{\nu,m}^{(p)}(x)\\
&+c\sum_{\nu>j} 2^{-|j-\nu|\varkappa}\left(M\left(\sum_{m} |\lambda_{\nu,m}| \chi_{\nu,m}^{(p)}\right)^w\right)^{\frac{1}{w}}(x).
\end{align*}
Now we can apply lemma $\ref{Delta}$ with 
\begin{align*}
 g_{\nu}:=\sum_{m}|\lambda_{\nu,m}| \tilde{\chi}_{\nu,m}^{(p)}(x)+\left(M\left(\sum_{m} |\lambda_{\nu,m}| \chi_{\nu,m}^{(p)}\right)^w\right)^{\frac{1}{w}}(x).
\end{align*}
Then it follows (with triangle inequality and the ``almost``-disjointness of the $c \cdot Q_{\nu,m}$)
\begin{align*}
 &\left\|\left(\sum_{j=0}^{\infty}2^{jsq}\left\|\left(k^N_j * \sum_{\nu,m} \lambda_{\nu,m} a_{\nu,m}\right)\big|E\right\|^q\right)^{\frac{1}{q}}\! \big|L_p\right\| \\
&\hspace{1em}\leq C \left\|\left(\sum_{\nu,m} \left(|\lambda_{\nu,m}|\ \tilde{\chi}_{\nu,m}^{(p)}\right)^q \right)^{\frac{1}{q}}\!\big|L_p\right\| 
+C\left\|\left(\sum_{\nu} \left(M\left(\sum_{m} |\lambda_{\nu,m}| \chi_{\nu,m}^{(p)}\right)^w\right)^{\frac{q}{w}}\right)^{\frac{1}{q}}\!\big|L_p\right\|
\end{align*}
and an analogous result for $B_{p,q}^s(E)$. Furthermore, the first term can be estimated by a term similar to the second, observing $\tilde{\chi}_{\nu,m}^{(p)}\leq C' \left(M\ \left(\chi_{\nu,m}^{(p)}\right)^w\right)^{\frac{1}{w}}$.

Eventually, the proposition follows by $\| \ \|M(f_n^w)^{\frac{1}{w}}|l_q\| \ |L_p\|=\| \ \|M(f_n^w)|l_{\frac{q}{w}}\| \ |L_{\frac{p}{w}}\|^{\frac{1}{w}}$, $w<p$, $w<q$ and by the boundedness of the maximal operator (see $\eqref{GrundMaxLplq}$).
In the case $B_{p,q}^s(E)$ we only need the boundedness of the maximal operator from $l_s(L_r)$ to $l_s(L_r)$, which is given for $r>1$ (see $\eqref{GrundMaxLp}$) such that $w<p$ and hence $L\geq \lfloor\sigma_p-s\rfloor$ suffices.
\end{bew}
Now the question will be whether all elements of the function space can be represented in such a way. The positive answer in the scalar (i.e. $E=\mathbb{C}$) case has been given for instance in \cite{Tri97}, theorem 13.8, p.$ $ 75. For the vector-valued case we will slightly alter the derivation sequence, as described in \cite{Tri97}, theorem 15.8, p.$ $ 114 . In the first step we care about a representation with harmonic, vector-valued atoms. This is inspired by the norms from proposition $\ref{RychkovHarmSatz}$, in which the functions $u(x,t)$ are harmonic in the domain $\left\{ x \in \mathbb{R}^n, t>0\right\}$. 

To explain this a bit more in detail (as in \cite{Tri97}, section 12.2, p.$ $ 59 in the scalar case) we choose an $f \in {\cal S}(\mathbb{R}^n,E)$ and form the functions $u(x,t)$ for $x \in \mathbb{R}^n$, $t>0$ as in example $\ref{RychkovHarmonisch}$ by
\begin{align*}
 u(x,t):=(e^{-t|\cdot|} \hat{f})\check{\ }(x)= d_n \left(f *  \frac{t}{(|\cdot|^2+t^2)^{\frac{n+1}{2}}}\right)(x).
\end{align*}
We obtain
\begin{align*}
 u(x,t) \rightarrow f(x) \text{ for } t \rightarrow 0
\end{align*}
uniformly in $x$ because $(e^{-|\cdot|})\check{\  } \in L_1$ and $e^{-|0|}=1$. 
Furthermore, we have
\begin{align*}
 t^k \frac{\partial^k u(x,t)}{\partial t^k} \rightarrow 0 \text{ for } t \rightarrow 0 
\end{align*}
uniformly in $x$ for all $k \in \mathbb{N}$ because $(|\cdot|^k e^{-|\cdot|})\check{\  } \in L_1$ and $|0|^ke^{-|0|}=0$.
Now we obtain by iterated partial integration and with suitable constants $d_l^k$ with $k \in \mathbb{N}$ and $l \in \{0,\ldots,k-1\}$
\begin{align*} 
\int_a^b t^{k-1} \frac{\partial^k u(x,t)}{\partial t^k} \ dt &= \tau^{k-1} \frac{\partial^{k-1} u(x,\tau)}{\partial \tau^{k-1}}\Big|_a^b - (k-1) \int_a^b t^{k-2} \frac{\partial^{k-1} u(x,t)}{\partial t^{k-1}} \\
&= \ldots \\
 &=\sum_{l=0}^{k-1} d_l^k b^l \frac{\partial^l u(x,b)}{\partial t^l}-\sum_{l=0}^{k-1} d_l^k a^l \frac{\partial^l u(x,a)}{\partial t^l}.
\end{align*}
Therefore, we get the relation
\begin{align*}
\sum_{l=0}^{k-1} d_l^k \tau^l \frac{\partial^l u(x,\tau)}{\partial t^l}\Big|_0^1 =\sum_{\nu=0}^{\infty}  \int_{2^{-\nu-1}}^{2^{-\nu}} t^{k-1} \frac{\partial^k u(x,t)}{\partial t^k} \ dt.
\end{align*}
By our above considerations on the limits it follows 
\begin{align}
\label{HarmDarst}
 f(x)=c\sum_{\nu=0}^{\infty} \int_{2^{-\nu-1}}^{2^{-\nu}} t^{k-1} \frac{\partial^k u(x,t)}{\partial t^k} \ dt+\sum_{l=0}^{k-1} c_l^k \frac{\partial^l u(x,1)}{\partial t^l}
\end{align}
with suitable constants $c_l^k$ with $k \in \mathbb{N}$ and $l \in \{0,\ldots,k-1\}$. We want to call the right-hand side a harmonic representation of $f$. A look at the norms from proposition $\ref{RychkovHarmSatz}$ tells us that very similar terms occured there. Therefore, it will be our aim to give $\eqref{HarmDarst}$ a meaning for $f \in B_{p,q}^s(E)$ resp.$ $ $F_{p,q}^s(E)$, with convergence at least in ${\cal S}'(\mathbb{R}^n,E)$. From the remarks forward to $\ref{RychkovHarmSatz}$ we obtain that the functions $u(x,t)$ are well-defined for $f \in B_{p,q}^s(E)$ resp.$ $ $F_{p,q}^s(E)$, harmonic in $\{(x,t) \in \mathbb{R}^{n+1}, t>0\}$ and bounded on $\{(x,t) \in \mathbb{R}^{n+1}, t>\delta\}$ for every $\delta>0$. So the integrals in $\eqref{HarmDarst}$ make sense. 

In the following we keep close to \cite{Tri97}, theorem 12.5, p. 62, where the scalar case is treated.
\begin{Satz}
\label{HarmDarstSa}
 Let $s \in \mathbb{R}$, $0<q\leq \infty$ and $0<p\leq \infty$ (resp. $\!<\infty$). If one chooses $k \in\mathbb{N}$ large enough, then the right-hand side of $\eqref{HarmDarst}$ converges in ${\cal S}'(\mathbb{R}^n,E)$ to $f$ for $f \in B_{p,q}^s(E)$ resp. $F_{p,q}^s(E)$.
\end{Satz}
\begin{bew}
This proposition can be proven as the scalar-case \cite{Tri97}, theorem 12.5(i), p.$ $ 62. To avoid (unknown) vector-valued duality relations one shows
\begin{align*}
 a[f]=a[g]
\end{align*}
for $a \in E'$ instead, where $a[f]$ with $a[f](\varphi):=a(f(\varphi))$ for $\varphi \in {\cal S}(\mathbb{R}^n)$ is an element of $B_{p,q}^s$ if $f \in B_{p,q}^s(E)$.
\end{bew}
Now we derive a representation of the elements of the function spaces $B_{p,q}^s(E)$ resp. $F_{p,q}^s(E)$ which leads us to the desired form of proposition 
$\ref{HarmIfAtom}$. Hereby we keep close to \cite{Tri97}, section 13.10, p.$ $ 83.

Let $f \in B_{p,q}^s(E)$ resp. $F_{p,q}^s(E)$. If we choose $k$ sufficiently large, then we obtain the harmonic representation of proposition $\ref{HarmDarstSa}$
\begin{align*}
 f(x)=c\sum_{\nu=0}^{\infty}  \int_{2^{-\nu-1}}^{2^{-\nu}} t^k \frac{\partial^k u(x,t)}{\partial t^k} \ \frac{dt}{t}+\sum_{l=0}^{k-1} c_l^k \frac{\partial^l u(x,1)}{\partial t^l}
\end{align*}
with convergence in ${\cal S}'(\mathbb{R}^n,E)$. 

Let $\mu \in \mathbb{N}$ be fixed. Let $\nu \in \mathbb{N}$, $\nu\geq\mu$, $m \in \mathbb{Z}^n$ and $l \in \{0,\cdots 2^{\mu}-1\}$. By $B_{\nu,m,l}$ we denote the cubes in $\mathbb{R}_+^{n+1}:=\{(x,t) \in \mathbb{R}^{n+1}, t>0 \}$ with center $(2^{-\nu}m,2^{-\nu+\mu}+l2^{-\nu})$ and radius $2^{-\nu+\mu-2}$. We decompose the rectangles $Q_{\nu,m} \times (2^{-\nu+\mu},2^{-\nu+\mu+1})$ in $2^{\mu}$ cubes of side length $2^{-\nu}$. Now we define
\begin{align}
\label{HarmLambda}
 \lambda_{\nu,m}:=2^{\nu\left(s-\frac{n}{p}\right)}2^{-\nu k} \sup \| \frac{\partial^k u(y,t)}{\partial t^k}|E\| \text{ for } \nu>\mu, m \in \mathbb{Z}^n,
\end{align}
where we take the supremum over the set
\begin{align*}
\left\{(y,t) \in \mathbb{R}^{n+1}: |2^{-\nu}m-y|\leq d2^{-\nu+\mu-1}, d^{-1} 2^{-\nu+\mu}\leq t \leq d 2^{-\nu+\mu+1}\right\},
\end{align*}
for a $d>0$ which we will choose sufficiently large afterwards. 

In the case $\mu=\nu$ we put
\begin{align*}
 \lambda_{\mu,m}:=\sum_{l=0}^{k-1}  |c_l^k|\sup \| \frac{\partial^l u(y,t)}{\partial t^l}|E\| \text{ for } m \in \mathbb{Z}^n,
\end{align*}
where we take the supremum over the set 
\begin{align*}
\left\{(y,t) \in \mathbb{R}^{n+1}: |2^{-\nu}m-y|\leq d, \frac{1}{2d} \leq t \leq \frac{3}{2}d \right\}.
\end{align*}
Now we take a closer look at $\|\lambda|f_{p,q}\|$. We obtain
\begin{align}
\label{HarmNormEq}
\begin{split}
 \|\lambda|f_{p,q}\|&=\left\|\left(\sum_{\nu=\mu}^{\infty} \sum_{m \in \mathbb{Z}^n} |\lambda_{\nu,m}\chi_{\nu,m}^{(p)}|^q\right)^{\frac{1}{q}}\big|L_p\right\| \\
&	\sim
 \left\|\left(\sum_{\nu=\mu+1}^{\infty} \sum_{m \in \mathbb{Z}^n} \left(2^{\nu\left(s-\frac{n}{p}\right)}2^{-\nu k} \sup \| \frac{\partial^k u(y,t)}{\partial t^k}|E\| \ \chi_{\nu,m}^{(p)}\right)^q\right)^{\frac{1}{q}}\big|L_p\right\| \\
&\hspace{2em}+\left\|\sum_{m \in \mathbb{Z}^n} \sum_{l=0}^{k-1}  |c_l^k|\sup \| \frac{\partial^l u(y,t)}{\partial t^l}|E\| \cdot\chi_{\mu,m}^{(p)} \ \big|L_p\right\| \\
&\leq c2^{\mu\delta} \|f|F_{p,q}^s(E)\|
\end{split}
\end{align}
for a suitable $\delta>0$ and for a $c$ independent of $\mu$. Note that the estimate of the norms by proposition $\ref{RychkovHarmLast}$ was applied in the last step of the chain of proof. Hence we have to assume $p>\frac{n}{n+1}$. 

Is it neglectable that the suprema in the second last term are taken over a small area of $t$ resp. a larger area of $t$: This is obvious for the second part and follows for the first part by a Taylor expansion. 

Analogously, $\|\lambda|b_{p,q}\|$ can be estimated by $c\|f|B_{p,q}^s(E)\|$ in the case $f \in B_{p,q}^s(E)$. Now we choose a $\psi \in {\cal S}(\mathbb{R}^n)$ with compact support and
\begin{align}
\label{HarmLocal}
 \sum_{m \in \mathbb{Z}^n} \psi(x-m)=1 \text{ for all } x \in \mathbb{R}^n.
\end{align}
If $\nu>\mu$ and $m\in \mathbb{Z}^n$, we put (with $c$ out of $\eqref{HarmDarst}$)
\begin{align*}
  a_{\nu,m}(x):=\! \sum_{l=0}^{2^{\mu}-1} a_{\nu,m,l}(x)\text{ with } \, a_{\nu,m,l}(x):=c \lambda_{\nu,m}^{-1} \psi(2^{\nu}x-m) \int_{2^{-\nu+\mu}+l2^{-\nu}}^{2^{-\nu+\mu}+(l+1)2^{-\nu}} \!t^k \frac{\partial^k u(x,t)}{\partial t^k} \ \frac{dt}{t}
\end{align*}
and in the case $\nu=\mu$ we define for $m \in \mathbb{Z}^n$ 
\begin{align*}
 a_{\mu,m}(x):=\sum_{l=0}^{k-1} a_{\mu,m,l}(x)\text{ with } \, a_{\mu,m,l}(x):=c_l^k \lambda_{\mu,m}^{-1} \psi(2^{\mu}x-m)\frac{\partial^l u(x,1)}{\partial t^l}.
\end{align*}
Then we obtain (in ${\cal S}'(\mathbb{R}^n,E)$)
\begin{align*}
 f&=\sum_{\nu=\mu}^{\infty} \sum_{m \in \mathbb{Z}^n} \lambda_{\nu,m} a_{\nu,m}.
\end{align*}
This is the desired representation. In the following we are going to show that the $a_{\nu,m}$ behave like $E$-valued $(s,p)_{K,-1}$-atoms for all $K \in \mathbb{N}$. 

By construction the condition $\eqref{Atom1}$ is valid. We can't show any moment conditions (see $\eqref{Atom3}$). To check the conditions on the derivatives we use a lemma for harmonic functions.
\begin{Lemma}
\label{HarmHarm}
 Let $W(X_1,\cdots,X_N): \mathbb{R}^N \rightarrow E$ be harmonic in the domain
\begin{align*}
 K_R=\left\{X \in \mathbb{R}^N: |X|\leq R\right\}.
\end{align*}
 Then for $\varkappa \in (0,1)$ it is true that
\begin{align*}
 \|D^{\alpha} W(X)|E\| \leq c_{\alpha,\varkappa} R^{-|\alpha|} \sup_{|Y|=R} \|W(Y)|E\| \text{ for } |X|\leq \varkappa R
\end{align*}
 with a constant $c$ which depends on $\alpha$ and $\varkappa$ but not on $R$.
\end{Lemma}
\begin{bew}
If $V: \mathbb{R}^N\rightarrow E$ is harmonic in the given domain, then it holds 
\begin{align}
\label{HarmDirichlet}
 V(X)=\frac{R^2-|X|^2}{R\omega_N} \int_{|Y|=R} \frac{V(Y)}{|X-Y|^N} \ ds_Y \text{ for } |X|<R,
\end{align}
where $\omega_N$ is the volume of the unit ball of $\mathbb{R}^N$. The lemma follows by taking the derivative of both sides and by a suitable estimate (see the end of \cite{Tri97}, section 13.10, p.$ $ 83).
\end{bew}
Now we apply this lemma to the functions $W(X)=\frac{\partial^k u(x,t)}{\partial t^k}$, which are harmonic in $\mathbb{R}^{n+1}_+$, to the balls $B_{\nu,m,l}$ instead of $K_R$ with $R=2^{-\nu+\mu-2}$ and to $\varkappa=d'2^{-\mu+2}<1$ if $\mu$ is larger than a certain $\varkappa_0$. Thus we obtain for $\nu>\mu$ and the set of all
\begin{align*}
\{(x,t) \in \mathbb{R}^{n+1}: |(2^{-\nu}m-x ,2^{-\nu+\mu}+l2^{-\nu}-t)|\leq d'2^{-\mu+2}\cdot 2^{-\nu+\mu-2}=d' 2^{-\nu}\}
\end{align*}
the relation
\begin{align*}
 \left\|D^{\gamma} \frac{\partial^k u(x,t)}{\partial t^k}|E\right\| &\leq c\, 2^{(\nu-\mu+2) |\gamma|} \sup\left\|\frac{\partial^k u(x,t)}{\partial t^k}|E\right\|,
\end{align*}
where the supremum is taken over
\begin{align*}
\{(x,t) \in \mathbb{R}^{n+1}: |(2^{-\nu}m-x ,2^{-\nu+\mu}+l2^{-\nu}-t)|=2^{-\nu+\mu-2}\}.
\end{align*}
But for all $l\in \{0,\ldots,2^{\mu}-1\}$ this set is contained in the set $(x,t) \in \mathbb{R}^{n+1}$ with $|2^{-\nu}m-x|\leq d2^{-\nu+\mu-1}$ and $d^{-1} 2^{-\nu+\mu}\leq t \leq d 2^{-\nu+\mu+1}$ for a suitable $d$. Now this yields
\begin{align*}
 \|D^{\gamma} a_{\nu,m}|E\|\leq c\sum_{l=0}^{2^{\mu}-1} 2^{-\nu\left(s-\frac{n}{p}\right)+\nu |\gamma|} 2^{\nu k} \int_{2^{-\nu+\mu}+l2^{-\nu}}^{2^{-\nu+\mu}+(l+1)2^{-\nu}} t^k \ \frac{dt}{t} 
\leq c' 2^{\mu k} 2^{-\nu(s-\frac{n}{p})+\nu |\gamma|}.
\end{align*}
Analogous assertions hold true for $a_{\mu,m}$ ($m \in \mathbb{Z}^n$). Therefore, we have proven the desired conditions $\eqref{Atom2}$ for all $K \in \mathbb{N}_0$. The $a_{\nu,m}$ introduced above for $\nu \in \mathbb{N}$, $\nu\geq\mu$ and $m \in \mathbb{Z}^n$ are $E$-valued $(s,p)_{K,-1}$-atoms for all $K \in \mathbb{N}_0$ - up to a constant depending on $\mu$.

We call this atoms and the found representation for $f$ ``harmonic``.
\begin{Satz}
\label{HarmAtomDarst}
 (i) Let $\frac{n}{n+1}<p\leq \infty$, $0<q\leq \infty$, $s>\sigma_p$ and $K \in \mathbb{N}_0$ with $K\geq 1+\lfloor s \rfloor$.
Then $f \in {\cal S}'(\mathbb{R}^n,E)$ belongs to $B_{p,q}^s(E)$ if and only if it can be represented by
\begin{align*}
 f=\sum_{\nu \in \mathbb{N}_0} \sum_{m \in \mathbb{Z}^n} \lambda_{\nu,m} a_{\nu,m}(x).
\end{align*}
Here $a_{\nu,m}$ are $E$-valued $1_K$-atoms (for $\nu=0$) or $E$-valued $(s,p)_{K,-1}$-atoms (for $\nu \in \mathbb{N}$) and $\lambda \in b_{p,q}$.
Furthermore, we have 
\begin{align*}
 \|f|B_{p,q}^s(E)\|  \sim \inf \|\lambda|b_{p,q}\|
\end{align*}
in the sense of equivalence of norms, where the infimum on the right-hand side is taken over all admissible representations for $f$.

 (ii) Let $\frac{n}{n+1}<p<\infty$, $0<q\leq \infty$, $s>\sigma_{p,q}$ and $K \in \mathbb{N}_0$ with $K\geq 1+\lfloor s \rfloor $.
Then $f \in {\cal S}'(\mathbb{R}^n,E)$ belongs to $F_{p,q}^s(E)$ if and only if it can be represented by
\begin{align*}
 f=\sum_{\nu \in \mathbb{N}_0} \sum_{m \in \mathbb{Z}^n} \lambda_{\nu,m} a_{\nu,m}(x).
\end{align*}
Here $a_{\nu,m}$ are $E$-valued $1_K$-atoms (for $\nu=0$) or $E$-valued $(s,p)_{K,-1}$-atoms (for $\nu \in \mathbb{N}$) and $\lambda \in f_{p,q}$.
Furthermore, we have 
\begin{align*}
 \|f|F_{p,q}^s(E)\|  \sim \inf \|\lambda|f_{p,q}\|
\end{align*}
in the sense of equivalence of norms, where the infimum on the right-hand side is taken over all admissible representations for $f$.
\end{Satz}
\begin{bew}
 The assertions follow from $\ref{HarmIfAtom}$ because the choice of $L=-1$ is admissible in the case $s>\sigma_p$ resp. $s>\sigma_{p,q}$ and from the proven representation observing $\eqref{HarmNormEq}$ resp.$ $ the analogous result for $B_{p,q}^s(E)$. Here the coefficients $\lambda_{\nu,m}$ even vanish for $\nu<\mu$.
\end{bew}

\subsection{Subatomic decompositions}
The aim of the following section will be to simplify the atomic representation of $f \in B_{p,q}^s(E)$ resp. $F_{p,q}^s(E)$ further. As a basis we use the harmonic representation from the last section. We orientate on \cite{Tri97}, section 14, p.$ $ 89 which treats the scalar case. 
\begin{Definition}
 Let $\psi \in {\cal S}(\mathbb{R}^n)$ with $supp \ \psi \subset d \cdot Q_{0,0}$ for a $d>1$ and 
\begin{align*}
 \sum_{m \in \mathbb{Z}^n} \psi(x-m)=1.
\end{align*}
Let $s \in \mathbb{R}$, $0<p\leq \infty$, $\frac{L+1}{2} \in \mathbb{N}_0$ and $\gamma \in \mathbb{N}_0^n$. We put $\psi^{\gamma}(x):=x^{\gamma} \psi(x)$. Then we call
\begin{align*}
(\gamma qu)_{\nu,m}^L(x)=2^{-\nu\left(s-\frac{n}{p}\right)} \left( (-\Delta)^{\frac{L+1}{2}} \psi^{\gamma}\right)(2^{\nu} x-m)
\end{align*}
 a $(s,p)_L$-$\gamma$-quark for $Q_{\nu,m}$. If $L=-1$, we want to denote it shortly by $(\gamma qu)_{\nu,m}(x)$. 
\end{Definition}
\begin{Bemerkung}
\label{HarmQuarksAtoms}
 First of all, we want to show that the $(s,p)_L-\gamma$-quarks really are (scalar) $(s,p)_{K,L}$-atoms for all $K \in \mathbb{N}_0$. The moment conditions $\eqref{Atom3}$ easily follow from their shape. For the derivatives we have
\begin{align*}
  \left|D^{\alpha} \left(\left((-\Delta)^{\frac{L+1}{2}} \psi^{\gamma}\right)(2^{\nu} x-m)\right)\right| &\leq c\, 2^{|\alpha|\nu}2^{\varkappa |\gamma|},
\end{align*}
where $c$ and $\varkappa$ depend on $\alpha$ and $L$ but not on $\gamma$, $\nu$, or $m$. So the $(\gamma qu)_{\nu,m}^L(x)$ are $(s,p)_{K,L}$-atoms up to a constant.
\end{Bemerkung}
Now we will simplify the representation of $f \in B_{p,q}^s(E)$ resp. $F_{p,q}^s(E)$ by the following result which corresponds to \cite{Tri97}, theorem 14.4, step 2, p.$ $ 93.
\begin{Lemma}
\label{HarmQuarkDarst}
(i) Let $\frac{n}{n+1}<p\leq \infty$, $0<q\leq\infty$, $s>\sigma_p$ and $f \in B_{p,q}^s(E)$. Then there is a $\varkappa_0 \in \mathbb{N}$ such that there exists a representation
\begin{align*}
 f=\sum_{\gamma \in \mathbb{N}_0^n} \sum_{\nu=0}^{\infty} \sum_{m \in \mathbb{Z}^n} \lambda_{\nu,m}^{\gamma} e_{\nu,m}^{\gamma} (\gamma qu)_{\nu,m}(x)
\end{align*}
in ${\cal S}'(\mathbb{R}^n,E)$ with $e_{\nu,m}^{\gamma} \in U_E$ for all $\mu \geq \varkappa_0$. It holds
\begin{align*}
 \sup_{\gamma \in \mathbb{N}_0^n} 2^{\mu|\gamma|} \|\lambda^{\gamma}|b_{p,q}\|\leq c \|f|B_{p,q}^s(E)\|,
\end{align*}
where $c$ does not depend on $f$ and $\lambda^{\gamma}=(\lambda_{\nu,m}^{\gamma})_{m \in \mathbb{Z}^n, \nu \in \mathbb{N}_0}$.

(ii) Let $\frac{n}{n+1}<p<\infty$, $0<q\leq\infty$, $s>\sigma_{p,q}$ and $f \in F_{p,q}^s(E)$. Then there is a $\varkappa_0 \in \mathbb{N}$ such that there exists a representation
\begin{align*}
 f=\sum_{\gamma \in \mathbb{N}_0^n} \sum_{\nu=0}^{\infty} \sum_{m \in \mathbb{Z}^n} \lambda_{\nu,m}^{\gamma} (\gamma qu)_{\nu,m}(x) e_{\nu,m}^{\gamma},
\end{align*}
in ${\cal S}'(\mathbb{R}^n,E)$ with $e_{\nu,m}^{\gamma} \in U_E$ for all $\mu \geq \varkappa_0$. It holds
\begin{align*}
 \sup_{\gamma \in \mathbb{N}_0^n} 2^{\mu|\gamma|} \|\lambda^{\gamma}|f_{p,q}\|\leq c \|f|F_{p,q}^s(E)\|,
\end{align*}
where $c$ does not depend on $f$ and $\lambda^{\gamma}=(\lambda_{\nu,m}^{\gamma})_{m \in \mathbb{Z}^n, \nu \in \mathbb{N}_0}$.
\end{Lemma}
\begin{bew}
 We restrict ourselves to the case $f \in F_{p,q}^s(E)$, the case $f \in B_{p,q}^s(E)$ can be proven analogously. From proposition $\ref{HarmAtomDarst}$ and the previous remarks  we obtain the optimal decomposition
\begin{align}
  \label{HarmDarst3}
 f=\sum_{\nu \geq \mu} \sum_{m \in \mathbb{Z}^n} \lambda_{\nu,m} a_{\nu,m}
\end{align}
with (in case of $\nu>\mu$)
\begin{align}
 \label{HarmAtomDef}
\begin{split}
a_{\nu,m,l}(x)&=c \lambda_{\nu,m}^{-1} \psi(2^{\nu}x-m) \int_{2^{-\nu+\mu}+l2^{-\nu}}^{2^{-\nu+\mu}+(l+1)2^{-\nu}} t^k \frac{\partial^k u(x,t)}{\partial t^k} \ \frac{dt}{t}, \\
  a_{\nu,m}(x)&= \sum_{l=0}^{2^{\mu}-1} a_{\nu,m,l}(x) \quad \text{and} \quad  \|\lambda|f_{p,q}\| \leq c2^{\mu \delta}  \|f|F_{p,q}^s(E)\|.
\end{split}
\end{align}
 We want to expand the arbitrarily often differentiable functions $\frac{\partial^k u(x,t)}{\partial t^k}$ into a Taylor series, with center $(2^{-\nu}m,2^{-\nu+\mu}+l2^{-v})$ of the balls $B_{\nu,m,l}$. We need 
\begin{Lemma}
\label{HarmHarmHarm}
 There exist $c>0$ and $0<\tau<1$ such that
\begin{align*}
  \|D^{\alpha} W(0)|E\| \leq c \alpha! \tau^{-|\alpha|} \sup_{|y|=1} \|W(y)|E\|
\end{align*}
for all $\alpha \in \mathbb{N}_0^N$ and all $W:\mathbb{R}^N \rightarrow E$ which are harmonic in the domain $\{y \in \mathbb{R}^N: |y|\leq 1\}$.
\end{Lemma}
\begin{bew}
For an arbitrary $Y \in \mathbb{R}^N$ with $|Y|=1$ we expand the function   
\begin{align*}
 |Z-Y|^{-N}=\left[ \sum_{j=1}^N (Z_j-Y_j)^2\right]^{-\frac{N}{2}},
\end{align*}
which is holomorphic in $|Z|<c$ with $c$ independent of $Y$, into its Taylor series around $0$. By a repeated application of the Cauchy formula for $\tau\leq c'<\frac{c}{\sqrt{N}}$ it follows
\begin{align*}
 \left|\left(D^{\alpha} |Z-Y|^{-N}\right)(0)\right|&=\left|\frac{(-1)^{\alpha}\alpha! }{(2\pi i)^{N}} \int_{|z_1|= \tau} \ldots \int_{|z_N|= \tau} \frac{|(z_1,\ldots,z_N)-Y|^{-N}}{z_1^{\alpha_1+1}\cdot\ldots\cdot z_N^{\alpha_N+1}} \ dz_1 \ldots dz_N \right|\\
&\leq c\tau^{-|\alpha|} \alpha!
\end{align*}
uniformly in $Y \in \mathbb{R}^N$ with $|Y|=1$.

By the formula of Dirichlet for $E$-valued functions which are harmonic in $\{Y \in \mathbb{R}^N:|Y|\leq 1\}$ (see $\eqref{HarmDirichlet}$) and by the uniform convergence of the Taylor series of $D^{\alpha} |Z-Y|^{-N}$ for $X,Y \in \mathbb{R}^N$ with $|Y|=1$ and $|X|<\tau$ we obtain
\begin{align*}
  W(X)&=\frac{1-|X|^2}{\omega_N} \int_{|Y|=1} \frac{W(Y)}{|X-Y|^N} \ ds_Y  \\
  &= \sum_{\alpha \in \mathbb{N}_0^N} \frac{X^{\alpha}}{\omega_N\alpha!} \int_{|Y|=1}  \tilde{a}_{\alpha}(Y) W(Y)  \ ds_Y
\end{align*}
with $a_{\alpha}(Y)=D^{\alpha} \left(|Z-Y|^{-N} \right) (0)$ and $\tilde{a}_{\alpha}(y)= a_{\alpha}(y)-\sum_{k=1}^N a_{\alpha-2e_k}(y)$
where $e_k$ is the multi-index $\left(0,\ldots,0,1,0,\ldots,0\right)$. By taking the derivative of the power series we see that
\begin{align*}
 \|(D^{\alpha}W)(0)|E\|&=\frac{1}{\omega_N} \left\|\int_{|Y|=1}  \tilde{a}_{\alpha}(Y) W(Y)  \ ds_Y \big|E\right\|
\leq c \alpha! \tau^{-|\alpha|} \sup_{|y|=1} \|W(y)|E\|
\end{align*}
is valid for all $\alpha \in \mathbb{N}_0^N$, where $c$ does not depend on $\alpha$ and $W$.
\end{bew}
Now we apply this lemma with $N=n+1$ to the functions $W(x,t):=\frac{\partial^k u(x,t)}{\partial t^k}$ which are harmonic in the ball $B_{\nu,m,l}$. If we set
\begin{align*}
\tilde{W}(x,t):= W\left(2^{-\nu+\mu-2}x+2^{-\nu}m,2^{-\nu+\mu-2}t+2^{-\nu+\mu}+l2^{-\nu}\right),
\end{align*}
then $\tilde{W}$ is harmonic in $\{y \in \mathbb{R}^{n+1}: |y|\leq 1 \}$. Hence we finally get the power series expansion
\begin{align}
\label{Harmhelp}
\begin{split}
 W(x,t)&=\tilde{W}(2^{\nu-\mu+2}x-2^{-\mu+2}m,2^{\nu-\mu+2}t-2^{2}-l2^{-\mu+2}) \\
&= \sum_{\alpha \in \mathbb{N}_0^N,\ \beta \in \mathbb{N}_0} c_{(\alpha,\beta)} 2^{(|\alpha|+\beta)(\nu-\mu+2)} \frac{(x-2^{-\nu}m)^{\alpha}\cdot(t-2^{-\nu+\mu}-l2^{-\nu})^{\beta}}{\alpha!\beta!}
\end{split}
\end{align}
for $|(x-2^{-\nu}m,t-2^{-\nu+\mu}-l2^{-\nu})|<2^{-\nu+\mu-2} \tau$.

If we choose $\mu$ larger or equal than a certain $\varkappa_0$, then this expansion is true in particular for $(x,t)\in \mathbb{R}^{n+1}$ with $x \in supp \ \psi(2^{\nu}x-m)$ and $t \in [2^{-\nu+\mu}+l2^{-\nu},2^{-\nu+\mu}+(l+1)2^{-\nu}]$.
Here we have
\begin{align*}
 \|c_{\alpha,\beta}|E\|=\|(D^{\alpha,\beta}\tilde{W}(0)|E\| \leq c \alpha! \tau^{-|\alpha|-\beta} \sup_{(x,t) \in B_{\nu,m,l}} \|W(y)|E\|
\end{align*}
by lemma $\ref{HarmHarmHarm}$ proven before. If we put this into $\eqref{HarmAtomDef}$, we obtain
\begin{align*}
 a_{\nu,m,l}(x)&=c \lambda_{\nu,m}^{-1} \psi(2^{\nu}x-m) \int_{2^{-\nu+\mu}+l2^{-\nu}}^{2^{-\nu+\mu}+(l+1)2^{-\nu}} t^k \frac{\partial^k u(x,t)}{\partial t^k} \ \frac{dt}{t} \\
&\equiv \sum_{\gamma \in \mathbb{N}_0^N} \tilde{c}_{\gamma} \ (2^{\nu}x-m)^{\gamma}\psi(2^{\nu}x-m)  
\end{align*}
with
\begin{align*}
  \|\tilde{c}_{\gamma}|E\|&\leq c 2^{-\nu|\gamma| } \lambda_{\nu,m}^{-1}\sum_{\beta=0}^{\infty} \frac{2^{(\nu-\mu+2)(|\gamma|+\beta)}}{\gamma!\beta!}\ c_{(\gamma,\beta)} \int\limits_{2^{-\nu+\mu}+l2^{-\nu}}^{2^{-\nu+\mu}+(l+1)2^{-\nu}} (t-2^{-\nu+\mu}-l2^{-\nu})^{\beta}t^k \ \frac{dt}{t} \\
&\leq c'2^{-\nu\left(s-\frac{n}{p}\right)} 2^{\mu k}  2^{(-\mu+2)|\gamma|} \tau^{-|\gamma|}\sum_{\beta=0}^{\infty} \tau^{-\beta} 2^{(-\mu+2) \beta},
\end{align*}
observing the definition of $\lambda_{\nu,m}$ in $\eqref{HarmLambda}$. The series over $\beta$ converges by our choice of $\mu$ and we get
\begin{align*}
  a_{\nu,m,l}(x)=\sum_{\gamma \in \mathbb{N}_0^n} \eta_{\nu,m,l}^{\gamma}(\gamma qu)_{\nu,m}(x) \, \, \text{ with }\,\, \| \eta_{\nu,m,l}^{\gamma}|E\|\leq c''2^{\mu k}(\tau^{-1} 2^{-\mu+2})^{|\gamma|}.
\end{align*}
Here $c''$ and $\tau$ are independent of $\mu$ and $\gamma$. If we replace $\mu$ by $M\mu$ afterwards, where $M \in \mathbb{N}$ is sufficiently large, and sum over $l=0,\ldots, 2^{\mu}-1$ in $\eqref{HarmAtomDef}$, we arrive at
\begin{align}
\label{HarmFixed}
a_{\nu,m}(x)=\sum_{\gamma \in \mathbb{N}_0^n} \eta_{\nu,m}^{\gamma}(\gamma qu)_{\nu,m}(x)\, \,
\text{ with } \, \, 
 \| \eta_{\nu,m}^{\gamma}|E\|\leq C2^{\mu \delta}2^{-\mu|\gamma|}
\end{align}
for certain $C>0$ and $\delta>0$ which do not depend on $\mu$ and $\gamma$. 

The case $\nu=\mu$ can be treated analogously. We just have to set $t=1$ in the Taylor expansion so that the sum over $\beta$ in $\eqref{Harmhelp}$ vanishes. 

Hence we obtain from $\eqref{HarmDarst3}$
\begin{align*}
 f=\sum_{\gamma \in \mathbb{N}_0^n} \sum_{\nu=0}^{\infty} \sum_{m \in \mathbb{Z}^n} \lambda_{\nu,m}^{\gamma} (\gamma qu)_{\nu,m}(x) e_{\nu,m}^{\gamma}
\end{align*}
in ${\cal S}'(\mathbb{R}^n,E)$ with 
\begin{align*}
\lambda_{\nu,m}^{\gamma}=\left\{\begin{array}{l l}
\lambda_{\nu,m} \|\eta_{\nu,m}^{\gamma}|E\|&,\nu\geq \mu \\
0&, \nu<\mu \end{array}
\right. \quad \text{ and } \quad
e_{\nu,m}^{\gamma}=
\left\{\begin{array}{l l}
\frac{\eta_{\nu,m}^{\gamma}}{\|\eta_{\nu,m}^{\gamma}|E\|}&,\nu\geq \mu \\
0&, \nu<\mu \end{array}
\right. .
\end{align*}
With $\eqref{HarmNormEq}$ and the observations on the dependence of $\mu$ in \eqref{HarmFixed} we find
\begin{align*}
 2^{\mu|\gamma|}\|\lambda^{\gamma}|f_{p,q}\|\leq C' 2^{\mu \delta_1}\|f|F_{p,q}^s(E)\|, \quad \gamma \in \mathbb{N}_0^n
\end{align*}
with $C'$ and $\delta_1$ independent of $\mu$ and $\gamma$ for $\mu \geq \varkappa_0$.
\end{bew}
Now we have all the ingredients together to prove \cite{Tri97}, theorem 15.8, p.$ $ 114, where now arbitrary $s \in \mathbb{R}$ are allowed.
\begin{Theorem}
 (i) Let $\frac{n}{n+1}<p\leq \infty$, $0<q\leq \infty$ and $s \in \mathbb{R}$. Let $M \in \mathbb{N}$ with $M>\sigma_p$ and $M>s$ and $L$ with $\frac{L+1}{2} \in \mathbb{N}_0$ and $L \geq \lfloor\sigma_p-s\rfloor$ be fixed. Let $(\gamma qu)_{\nu,m}$ and $(\gamma qu)_{\nu,m}^L$ be given as $(M,p)_{-1}$- resp. $(s,p)_{L}$-$\gamma$-quarks for a given function $\psi \in {\cal S}(\mathbb{R}^n)$ with compact support and the property $\eqref{HarmLocal}$. Then there exists a $\varkappa>0$ such that for all $\mu \geq \varkappa$ it is valid that $f \in {\cal S}'(\mathbb{R}^n,E)$ belongs to $B_{p,q}^s(E)$ if and only if it can be represented as
\begin{align*}
 f=\sum_{\gamma \in \mathbb{N}_0^n} \sum_{\nu=0}^{\infty} \sum_{m \in \mathbb{Z}^n}\varrho_{\nu,m}^{\gamma} e_{\nu,m}^{\gamma} (\gamma qu)_{\nu,m}(x) +\lambda_{\nu,m}^{\gamma} e_{\nu,m}^{\gamma,L} (\gamma qu)_{\nu,m}^L(x)
\end{align*}
in ${\cal S}'(\mathbb{R}^n,E)$ with $e_{\nu,m}^{\gamma},e_{\nu,m}^{\gamma,L} \in U_E$ and 
\begin{align*}
 \sup_{\gamma \in \mathbb{N}_0} 2^{\mu |\gamma|}(\|\varrho^{\gamma}|b_{p,q}\|+\|\lambda^{\gamma}|b_{p,q}\|) < \infty.
\end{align*}
Furthermore, it holds in the sense of equivalence of norms
\begin{align*}
 \|f|B_{p,q}^s(E)\| \sim \inf \ \sup_{\gamma} 2^{\mu |\gamma|}(\|\varrho^{\gamma}|b_{p,q}\|+\|\lambda^{\gamma}|b_{p,q}\|),
\end{align*}
where the inf. on the right-hand side is taken over all admissible representations of $f$.

(ii) Let $\frac{n}{n+1}<p<\infty$, $0<q\leq \infty$ and $s \in \mathbb{R}$. Let $M \in \mathbb{N}$ with $M>\sigma_{p,q}$ and $M>s$ and $L$ with $\frac{L+1}{2} \in \mathbb{N}_0$ and $L \geq \lfloor\sigma_{p,q}-s\rfloor$ be given. The quarks have the same meaning as in (i). Then there exists a $\varkappa>0$ such that for all $\mu \geq \varkappa$ it is valid that $f \in {\cal S}'(\mathbb{R}^n,E)$ belongs to $F_{p,q}^s(E)$ if and only if it can be represented as
\begin{align}
\label{HarmDarst4}
 f=\sum_{\gamma \in \mathbb{N}_0^n} \sum_{\nu=0}^{\infty} \sum_{m \in \mathbb{Z}^n}\varrho_{\nu,m}^{\gamma} e_{\nu,m}^{\gamma} (\gamma qu)_{\nu,m}(x)+ \lambda_{\nu,m}^{\gamma} e_{\nu,m}^{\gamma,L} (\gamma qu)_{\nu,m}^L(x)
\end{align}
in ${\cal S}'(\mathbb{R}^n,E)$ with $e_{\nu,m}^{\gamma},e_{\nu,m}^{\gamma,L} \in U_E$ and 
\begin{align}
\label{HarmNorm2}
 \sup_{\gamma} 2^{\mu |\gamma|}(\|\varrho^{\gamma}|f_{p,q}\|+\|\lambda^{\gamma}|f_{p,q}\|) < \infty.
\end{align} 
Furthermore, it holds in the sense of equivalence of norms
\begin{align*}
 \|f|F_{p,q}^s(E)\| \sim \inf \ \sup_{\gamma \in \mathbb{N}_0} 2^{\mu |\gamma|}(\|\varrho^{\gamma}|f_{p,q}\|+\|\lambda^{\gamma}|f_{p,q}\|),
\end{align*}
where the inf. on the right-hand side is taken over all admissible representations of $f$.
\end{Theorem}
\begin{bew}
We only consider the case $F_{p,q}^s(E)$. The proof for $B_{p,q}^s(E)$ can be organized analogously. Let $f \in {\cal S}'(\mathbb{R}^n,E)$ be represented by $\eqref{HarmDarst4}$ with the condition $\eqref{HarmNorm2}$. Then
\begin{align*}
 f_1^{\gamma}:=\sum_{\nu=0}^{\infty} \sum_{m \in \mathbb{Z}^n}\varrho_{\nu,m}^{\gamma} e_{\nu,m}^{\gamma} (\gamma qu)_{\nu,m}(x)
\end{align*}
and 
\begin{align*}
 f_2^{\gamma}:=\sum_{\nu=0}^{\infty} \sum_{m \in \mathbb{Z}^n} \lambda_{\nu,m}^{\gamma} e_{\nu,m}^{\gamma,L} (\gamma qu)_{\nu,m}^L(x)
\end{align*}
are represented as sums of atoms (up to a constant) as elements of ${\cal S}'(\mathbb{R}^n,E)$ by remark $\ref{HarmQuarksAtoms}$. Here $(\gamma qu)_{\nu,m}e_{\nu,m}^{\gamma}$ resp. $(\gamma qu)_{\nu,m}^Le_{\nu,m}^{\gamma,L}$ are $E$-valued $(M,p)_{K,-1}$- resp. $(s,p)_{K,L}$-atoms for every $K \in \mathbb{N}_0$, where one has to keep in mind a normalization constant $c 2^{\varkappa\gamma}$ with $c$ depending on $K$ but independent of $\gamma$ (see remark $\ref{HarmQuarksAtoms}$). Thus we obtain by proposition $\ref{HarmIfAtom}$ that $f_1^{\gamma} \in F_{p,q}^M(E)$\footnote{Because of $M>\sigma_{p,q}$ we need no moment conditions $\eqref{Atom3}$ for these atoms.}, that $f_2^{\gamma} \in F_{p,q}^s(E)$ and that there exists a $c''>0$ such that it holds
\begin{align*}
 \|f^{\gamma}|F_{p,q}^s(E)\|\leq c'
 \left(\|f_1^{\gamma}|F_{p,q}^M(E)\|+
 \|f_2^{\gamma}|F_{p,q}^s(E)\|\right)\leq c''2^{\varkappa|\gamma|} \left(\|\varrho^{\gamma}|f_{p,q}\| +\|\lambda^{\gamma}|f_{p,q}\|\right).
\end{align*}
with $f^{\gamma}=f_1^{\gamma}+f_2^{\gamma}$. Here $c''$ and $\varkappa$ are independent of $\gamma$ (and $f^{\gamma}$). Therefore, if we take $\mu>\varkappa$ for granted, it results from $\eqref{HarmNorm2}$ and a typical Minkowski/Hölder argument that
\begin{align*}
 f=\sum_{\gamma \in \mathbb{N}_0^n} f^{\gamma} \quad \text{ with } \quad  \|f|F_{p,q}^s(E)\|& \leq C \sup_{\gamma} 2^{\mu|\gamma|} \left(\|\varrho^{\gamma}|f_{p,q}\|+\|\lambda^{\gamma}|f_{p,q}\|\right)
\end{align*}
in $F_{p,q}^s(E)$. Hence this part of the proof is even valid for all $0<p\leq \infty$, $0<q\leq \infty$ and $s \in \mathbb{R}$.

Let $f$ from $F_{p,q}^s(E)$ be given for the second part of the proof. In the case $s>\sigma_{p,q}$ and $L=-1$ the assertion of the proposition follows from lemma $\ref{HarmQuarkDarst}$. Here we don't need any terms of the form $\varrho_{\nu,m}^{\gamma} e_{\nu,m}^{\gamma} (\gamma qu)_{\nu,m}(x)$.

Let now $s$ be arbitrary and $f \in F_{p,q}^s(E)$. Then we have by the lift property (see $\eqref{GrundLift}$)
\begin{align*}
 g=\left((1+|\cdot|^2)^{-\frac{L+1}{2}}\hat{f}\right)\check{\ } \in F_{p,q}^{s+L+1}(E)
\end{align*}
with $\|f|F_{p,q}^s(E)\| \sim \|g|F_{p,q}^{s+L+1}(E)\|$. Thus $f$ can be represented as
\begin{align*}
 f=g+(-\Delta)^{\frac{L+1}{2}} g.
\end{align*}
If we apply the same argument to $g$ and iterate the procedure, we obtain
\begin{align*}
 f=f_1+(-\Delta)^{\frac{L+1}{2}} f_2
\end{align*}
with $\|f|F_{p,q}^{s}(E)\| \sim \|f_1|F_{p,q}^{s+m(L+1)}(E)\|+\|f_2|F_{p,q}^{s+L+1}(E)\|$. If $L\geq \lfloor \sigma_{p,q}-s \rfloor$ (and $\frac{L+1}{2} \in \mathbb{N}_0$), then $s+L+1$ fulfils the conditions from lemma $\ref{HarmQuarkDarst}$, this means $s+L+1>\sigma_{p,q}$. Then $f_2$ can be represented by
\begin{align*}
 f_2=\sum_{\gamma \in \mathbb{N}_0^n} \sum_{\nu=0}^{\infty} \sum_{m \in \mathbb{Z}^n} \lambda_{\nu,m}^{\gamma} e_{\nu,m}^{\gamma,L} (\gamma qu)_{\nu,m}(x),
\end{align*}
where $(\gamma qu)_{\nu,m}$ are $(s+L+1,p)_{-1}$-$\gamma$-quarks and it holds
\begin{align*}
\|f_2|F_{p,q}^{s+L+1}(E)\| \sim \sup_{\gamma} 2^{\mu |\gamma|}\|\lambda^{\gamma}|f_{p,q}\|.
\end{align*}
But now we have
\begin{align*}
 (-\Delta)^{\frac{L+1}{2}}(\gamma qu)_{\nu,m}(x)&=(-\Delta)^{\frac{L+1}{2}} \left(2^{-\nu\left(s+L+1-\frac{n}{p}\right)} \psi^{\gamma} (2^{\nu} x-m) \right)\\
&= 2^{-\nu\left(s-\frac{n}{p}\right)} \left((-\Delta)^{-\frac{L+1}{2}} \psi^{\gamma} \right)(2^{\nu} x-m),
\end{align*}
which is an $(s,p)_L$-$\gamma$-quark. 

Furthermore, let's choose $m$ so large that $\tilde{M}:=s+m(L+1)$ fulfils the condition $\tilde{M}\geq M$. From $f_1 \in F_{p,q}^{\tilde{M}}(E)$ follows $f_1 \in F_{p,q}^{M}(E)$ as well. This yields a representation for $f_1$ with $(M,p)_{-1}$-$\gamma$-quarks by lemma $\ref{HarmQuarkDarst}$, observing $M>\sigma_{p,q}$. Hence we obtain a representation for $f$ as a sum of $(M,p)_{-1}$- and $(s,p)_L$-$\gamma$-quarks by
\begin{align*}
 f=\sum_{\gamma \in \mathbb{N}_0^n} \sum_{\nu=0}^{\infty} \sum_{m \in \mathbb{Z}^n}\varrho_{\nu,m}^{\gamma} e_{\nu,m}^{\gamma} (\gamma qu)_{\nu,m}(x)+\lambda_{\nu,m}^{\gamma} e_{\nu,m}^{\gamma,L} (\gamma qu)_{\nu,m}^L(x)
\end{align*}
and it holds by the previous steps 
\begin{align*}
\sup_{\gamma} 2^{\mu|\gamma|} \left(\|\varrho^{\gamma}|f_{p,q}\|+\|\lambda^{\gamma}|f_{p,q}\|\right)&\leq c \left(\|f_1|F_{p,q}^{M}(E)\|+\|f_2|F_{p,q}^{s+L+1}(E)\|\right)  \\
&\leq c'\|f|F_{p,q}^s(E)\|.
\end{align*}
\end{bew}
Now it is an easy task to expand this theorem to the more general atoms. This was suggested by the first step of the preceding proof in which we only used that the quarks are atoms. We now obtain \cite{Tri97}, theorem 15.11, p. 116.
\begin{Theorem}
 (i) Let $\frac{n}{n+1}<p\leq \infty$, $0<q\leq \infty$ and $s \in \mathbb{R}$. Let $M \in \mathbb{N}$ with $M>\sigma_p$ and $M>s$, $K \in \mathbb{N}_0$ with $K\geq \lfloor s \rfloor +1$ and $L$ with $\frac{L+1}{2} \in \mathbb{N}_0$ and $L \geq \lfloor\sigma_p-s\rfloor$ be fixed. Then there exists a $\varkappa>0$ such that for all $\mu\geq \varkappa$ it is valid that $f \in {\cal S}'(\mathbb{R}^n,E)$ belongs to $B_{p,q}^s(E)$ if and only if it can be represented by
\begin{align*}
 f=\sum_{\gamma \in \mathbb{N}_0^n} \sum_{\nu=0}^{\infty} \sum_{m \in \mathbb{Z}^n}\varrho_{\nu,m}^{\gamma} e_{\nu,m}^{\gamma} a_{\nu,m}^{\gamma}(x) +\lambda_{\nu,m}^{\gamma} e_{\nu,m}^{\gamma,L} a_{\nu,m}^{\gamma,L}(x)
\end{align*}
in ${\cal S}'(\mathbb{R}^n,E)$. Here $a_{\nu,m}^{\gamma}$ resp.$ $ $a_{\nu,m}^{\gamma,L}$ are $(M,p)_{K,-1}$ resp. $(s,p)_{K,L}$-atoms, $e_{\nu,m}^{\gamma},e_{\nu,m}^{\gamma,L} \in U_E$ and 
\begin{align*}
 \sup_{\gamma} 2^{\mu |\gamma|}(\|\varrho^{\gamma}|b_{p,q}\|+\|\lambda^{\gamma}|b_{p,q}\|) < \infty.
\end{align*}
Furthermore, we have in the sense of equivalence of norms
\begin{align*}
 \|f|B_{p,q}^s(E)\| \sim \inf \ \sup_{\gamma} 2^{\mu |\gamma|}(\|\varrho^{\gamma}|b_{p,q}\|+\|\lambda^{\gamma}|b_{p,q}\|),
\end{align*}
where the inf. on the right-hand side is taken over all admissible representations of $f$.

(ii) Let $\frac{n}{n+1}<p<\infty$, $0<q\leq \infty$ and $s \in \mathbb{R}$. Let $M \in \mathbb{N}$ with $M>\sigma_{p,q}$ and $M>s$, $K \in \mathbb{N}_0$ with $K\geq \lfloor s \rfloor +1$ and $L$ with $\frac{L+1}{2} \in \mathbb{N}_0$ and $L \geq \lfloor\sigma_{p,q}-s\rfloor$ be fixed. Then there exists a $\varkappa>0$ such that for all $\mu\geq \varkappa$ it is valid that $f \in {\cal S}'(\mathbb{R}^n,E)$ belongs to $F_{p,q}^s(E)$ if and only if it can be represented by
\begin{align*}
 f=\sum_{\gamma \in \mathbb{N}_0^n} \sum_{\nu=0}^{\infty} \sum_{m \in \mathbb{Z}^n}\varrho_{\nu,m}^{\gamma} e_{\nu,m}^{\gamma} a_{\nu,m}^{\gamma}(x) +\lambda_{\nu,m}^{\gamma} e_{\nu,m}^{\gamma,L} a_{\nu,m}^{\gamma,L}(x)
\end{align*}
in ${\cal S}'(\mathbb{R}^n,E)$. Here $a_{\nu,m}^{\gamma}$ resp. $a_{\nu,m}^{\gamma,L}$ are $(M,p)_{K,-1}$ resp.$ $ $(s,p)_{K,L}$-atoms, $e_{\nu,m}^{\gamma},e_{\nu,m}^{\gamma,L} \in U_E$ and 
\begin{align*}
 \sup_{\gamma} 2^{\mu |\gamma|}(\|\varrho^{\gamma}|f_{p,q}\|+\|\lambda^{\gamma}|f_{p,q}\|) < \infty.
\end{align*}
Furthermore, we have in the sense of equivalence of norms
\begin{align*}
 \|f|F_{p,q}^s(E)\| \sim \inf \ \sup_{\gamma} 2^{\mu |\gamma|}(\|\varrho^{\gamma}|f_{p,q}\|+\|\lambda^{\gamma}|f_{p,q}\|),
\end{align*}
where the inf. on the right-hand side is taken over all admissible representations of $f$.
\end{Theorem}
\begin{bew}
  The existence of such a representation for $f \in B_{p,q}^s(E)$ resp. $f \in F_{p,q}^s(E)$ follows from the fact that the $(s,p)_L$-$\gamma$-quarks are also $(s,p)_{K,L}$-atoms for all $K \in \mathbb{N}_0$ and by the previous theorem.  If $f \in {\cal S}(\mathbb{R}^n,E)$ can be represented in the given way, then it belongs to $B_{p,q}^s(E)$ resp. $F_{p,q}^s(E)$ by the first step of the proof of the previous theorem because it only uses that the quarks are $(M,p)_{K,-1}$- resp. $(s,p)_{K,L}$-atoms.
\end{bew}

	\setbibpreamble{The bibliographical references are sorted by the names of the (first) authors. \par\bigskip}

\noindent Benjamin Scharf \\
E-Mail: benjamin.scharf@uni-jena.de \\
\\
Friedrich-Schiller-Universit\"at Jena \\
Fakult\"at f\"ur Mathematik und Informatik \\
Mathematisches Institut \\
D-07737 Jena 
	
\end{document}